\documentclass[A4,times,review,sort&compress]{elsarticle}
\usepackage[active]{srcltx}
\usepackage{bbm}
\usepackage{mathrsfs}
\usepackage{tipa}
\usepackage{amsthm}
\usepackage{amssymb}
\usepackage{stackrel}
\usepackage{amsmath,amscd}
\usepackage[numbers]{natbib}
\newtheorem{dl}{Theorem}[section]
\newtheorem{tl}[dl]{Corollary}
\newtheorem{yl}[dl]{Lemma}
\newtheorem{dy}[dl]{Definition}
\newtheorem{llz}[dl]{Example}
\newtheorem{lz}[dl]{Formula}

\newtheorem{remark}[dl]{Remark}

\numberwithin{equation}{section}
\usepackage{hyperref}
\numberwithin{equation}{section}
\newproof{pot1}{Proof of Theorem \ref{mainthnew}}
\newproof{pot2}{Proof of Theorem \ref{mainthnew2}}
\newproof{pot3}{Proof of Theorem \ref{qanalogue}}
\newproof{pot4}{Proof of Theorem \ref{333}}
\newcommand{\poq}[2]{(#1;q)_{#2}}

\newcommand{\poqq}[2]{(#1;q^2)_{#2}}

\newcommand{\be}{\begin{equation}}
\newcommand{\ee}{\end{equation}}
\newcommand{\ba}{\begin{array}}
\newcommand{\ea}{\end{array}}
\newcommand{\bmn}{\begin{eqnarray}}
\newcommand{\emn}{\end{eqnarray}}
\newcommand{\bnm}{\begin{eqnarray*}}
\newcommand{\enm}{\end{eqnarray*}}
\newcommand{\bln}{\begin{subequations}}
\newcommand{\eln}{\end{subequations}}

\begin{document}
%\begin{frontmatter}
\title{Expansion Formulas of  Basic Hypergeometric Series via  the \mbox{$(1-xy,y-x)$}-Inversion  and Their Applications}
\author{Jin Wang\fnref{fn3,fn4}}
\fntext[fn3]{This  work was supported by NSF of Zhejiang Province (Grant~No.~LQ20A010004) and by NSF of China (Grant~No.~12001492).}
\fntext[fn4]{E-mail address: jinwang@zjnu.edu.cn}
\address[Canada]{College of Mathematics and Computer Science, Zhejiang Normal University, Jinhua 321004, P.~R. ~China}
\author{Xinrong Ma\fnref{fn1,fn2}}
\fntext[fn1]{This  work was supported by NSFC grant No. 11971341.}
\fntext[fn2]{E-mail address: xrma@suda.edu.cn.}
\address[P.R.China]{Department of Mathematics, Soochow University, Suzhou 215006, P.R.China}

\begin{abstract}
With the use of  the  $(f,g)$-matrix inversion under specializations that $f=1-xy,g=y-x$, we establish an $(1-xy,y-x)$-expansion formula. When specialized to basic hypergeometric series, this $(1-xy,y-x)$-expansion formula leads us to some expansion formulas expressing any  ${}_{r}\phi_{s}$ series in variable $x~t$ in terms of a linear combination of  ${}_{r+2}\phi_{s+1}$ series in $t$, as well as various specifications.  All these results can be regarded as common generalizations of many known expansion formulas in the setting of $q$-series.  As specific applications, some new transformation formulas of $q$-series including  new approach to the Askey-Wilson polynomials, the Rogers-Fine identity,  Andrews' four-parametric  reciprocity theorem and Ramanujan's ${}_1\psi_1$ summation formula, as well as a transformation for certain well-poised Bailey pairs, are presented.
\end{abstract}
\begin{keyword} $(f,g)$-matrix inversion, $(1-xy,y-x)$-expansion formula, basic hypergeometric series, transformation, summation, WP Bailey pair, Rogers-Fine identity, reciprocity theorem,  Ismail's argument,   $q$-series.

\vspace{15pt}
\end{keyword}
%%%%%%%%%%%%%%%%%%%%%%%%%%%%%%%%%%%%%%%%%%%%%%%%%%%%%%%%%%%%%%%%%%%%%%%
\maketitle
%%%%%%%%%%%%%%%%%%%%%%%%%%%%%%%%%%%%%%%%%%%%%%%%%%%%%%%%%%%%%%
%\tableofcontents
%%%%%%%%%%%%%%%%%%%%%%%%%%%%%%%%%%%%%%%%%%%%%%%%%%%%%%%%%%%%%%
%\vspace{20pt}
%\parskip 7pt
%\baselineskip 16pt
\section{Introduction and main results}

It is well known that the core of the classical Lagrange inversion
formula  \cite[cf.][Appendix E]{andrews4}  and \cite[cf.][\S 7.32]{1001} is to express  the coefficients
$a_n$ in the expansion of
\begin{eqnarray}
  F(x)=\sum_{n=0}^{\infty}a_n\left(\frac{x}{\phi(x)}\right)^n,\label{lag}
\end{eqnarray}
where $a_0=F(0)$ and for $n\geq 1$,
$$
a_n=\frac{1}{n!}\frac{d^{n-1}}{dx^{n-1}}\left[\phi^{n}(x)\frac{dF(x)}{dx}\right]_{x=0}
$$
provided that $F(x)$ and $\phi(x)$ are analytic around $x=0$, $
\phi(0)\neq 0$, $\frac{d}{dx}$ denotes the usual derivative
operator.

 In the past decades, $q$-analogues (as generalizations) of the Lagrange inversion formula have drawn  a lot of attentions
 \cite[cf.][]{andrew3,6,111,kratt,18,22,11}. For a
good survey about
 results and open problems on this topic, we prefer to refer the
reader to Stanton's paper \cite{11} and  only record here, for
 comparison purpose, four relevant highlights.

The first one is a $q$-analogue  found by Carlitz
\cite[cf.][Eq.(1.11)]{carlitz} in 1973, subsequently reproduced by Roman
\cite[cf.][p.253, Eq.(8.4)]{3rr} via $q$-umbral calculus, asserting that
for any formal
 series $F(x)$, it holds
\begin{eqnarray}
F(x)=\sum_{n=0}^{\infty}\frac{x^n}{(q,x;q)_n}\left[D^n_{q,x}\{F(x)(x;q)_{n-1}\}\right]_{x=0},\label{carbe}
\end{eqnarray}
where  $D_{q,x}$ denotes the usual $q$-derivative operator defined by
\[D_{q,x}\{F(x)\}=\frac{F(x)-F(xq)}{x(1-q)}.\]

Another path-breaking result is certainly attributed to  Gessel and Stanton. With the mind that
 the essence of the Lagrange inversion
formula is equivalent to finding  a pair of matrix inversions,
 they successfully discovered   many $q$-analogues of  the Lagrange inversion formula in \cite{111}, e.g., Theorems 3.7 and 3.15  therein.  We pause at this point to recall the concept of matrix inversions. Following  \cite{18,milne,schlosser2}, a matrix inversion is usually  defined to be a pair of infinite lower
triangular matrices $F=(A_{n,k})_{n,k\in {\mathbb N}}$  and $G=(B_{n,k})_{n,k\in {\mathbb N}}$, ${\mathbb N}$ is the set of nonnegative integers, such that  $A_{n,k}=B_{n,k}=0$ unless $\
n\geq k$ and $A_{n,n}B_{n,n}\neq 0$, and
$$\sum_{n\geq i\geq
k}A_{n,i}B_{i,k}=\sum_{n\geq i\geq
k}B_{n,i}A_{i,k}=\delta_{n,k},
$$
where $\delta_{n,k}$ denotes the usual Kronecker delta. From this perspective,  Gessel and Stanton established the following $q$-Lagrange inversion formula:
\begin{align} F(x)=\sum_{n\geq k\geq
0} a_k\frac{(Ap^kq^k;p)_{n-k}}{\poq{q}{n-k}}q^{-nk}x^n\label{gessel1}
\end{align}
if and only if
\begin{align}
a_n=\sum_{k=0}^n(-1)^{n-k}q^{\binom{n-k+1}{2}+nk}
\frac{(1-Ap^kq^k)(Aq^np^{n-1};p^{-1})_{n-k-1}}{\poq{q}{n-k}}F(q^k).\label{gessel2}
\end{align}

The third important result is the following $q$-expansion formula  discovered in 2002
by Liu \cite[Theorem 2]{liu}. By the use of  the  $q$-derivative  operator, together with Carlitz's $q$-analogue   (\ref{carbe})  and  the very-well-poised   $\,_6\phi_5$ summation formula \cite[(II.20)]{10}, Liu showed that
\begin{eqnarray}
F(x)=\sum_{n=0}^{\infty}\frac{(1-aq^{2n})(aq/x;q)_nx^n}{(q;q)_n(x;q)_n}\left[D^n_{q,x}\{F(x)(x;q)_{n-1}\}
\right]_{x=aq} \label{1111}.
\end{eqnarray}

Likewise Liu but in a more systematical way,  Chu in his work \cite{chu1} investigated  various functions whose $n$th $q$-derivative   can be given in closed form, thereby leading to $q$-series identities. His main result is as follows:
\begin{eqnarray}
F(x)=\sum_{n=0}^{\infty}\frac{(1-abq^{2n+\epsilon})(a/x;q)_nx^n}{(q;q)_n}\frac{\left[D^n_{q,x}\{F(x)(bx;q)_{n+\epsilon}\}\right]_{x=a}}{\poq{bx}{n+1+\epsilon}} \label{2222}\quad(\epsilon=\pm 1).\label{chuexp}
\end{eqnarray}

As is known today,  these $q$-expansion formulas have been proved to be very
important to the theory of  basic hypergeometric series. Nevertheless,  as has been pointed out by Gasper in \cite{188} ``$\cdots$ the succeeding higher order derivatives
becomes more and more difficult to calculate for $|z|<1$, and so one
is forced to abandon this approach and to search for another way
$\cdots$". Truly, a
comparison of the aforementioned results shows   that the expansion (\ref{gessel1})/(\ref{gessel2}) of  Gessel and Stanton via the method of matrix inversions  has an advantage over other three expansions in that it  avoids  calculation of higher order $q$-derivatives. As such, is it possible to deal with the other three expansions  by the same way?

Another motivation comes from our observation that in their authoritative  book \cite[\S 2.2]{10}, Gasper and Rahman recorded a  general expansion formula expressing a terminating ${}_{r+4}\phi_{r+3}$ series as a finite linear combination of other terminating  ${}_{r+2}\phi_{r+1}$ series, i.e., Eq.(2.2.4) of \cite{10}. As a $q$-analogue of Bailey's formula  \cite[4.3(1)]{bailey}, this finite expansion serves as a ladder connecting a few  basic but most useful summation formulas of $q$-series. Similar results on hypergeometric series can be found in \cite{field} by Fields  and Wimp. We also have been attracted by some similar expansion formulas related to  the famous Askey-Wilson polynomials were presented in the very recent papers \cite[cf.][Theorem 2.2]{ismail3} by Ismail and Stanton, \cite[cf.][Theorem 1.5, Propostions 1.8 and 4.1]{liunew} by Liu, and \cite[cf.][Theorem 2.6]{jia-zeng} by Jia and Zeng.
Therefore, it comes no surprise to ask whether any ${}_{r}\phi_{s}$ series be expressed, by reversing the direction of the above ladder if necessary, as a linear combination of higher   ${}_{r+m}\phi_{s+n}$ series (integers $m,n\geq 1$).  Some supporting evidence for this question is expansion formulas stated   as Theorems 2.8 and 2.9 in \cite{ismail3}.

To answer these two problems in full generality, ultimately finding out new $q$-analogues of the Lagrange inversion formula \cite[cf.][]{111, 11},  forms  the main theme of the present  paper. One of valid tools for this purpose, in our viewpoint, is  the so-called  $(f,g)$-matrix inversion initially appeared
in \cite{0020}, together with an important method often referred to  as  ``Ismail's argument" \cite{schlosser1}
in the literature.  So named because it is Ismail who was the first to show  Ramanujan's $\,_1\psi_1$ summation formula
\cite{ismail0} by analytic continuation \cite{ahl}. Later,  in \cite{ask} Askey and Ismail  used the same method  and
 Rogers' very-well-poised  $\,_6\phi_5$ summation formula \cite[(II.20)]{10}  to evaluate
Bailey's bilateral  $\,_6\psi_6$ series. The reader may consult  \cite{kadell} by Kadell for a more systematic exposition on applications of this method in the theory of $q$-series.  For completeness, we rephrase  Ismail's argument   by the following uniqueness theorem of analytic functions \cite[cf.][p.127]{ahl}.

\begin{yl}[Uniqueness theorem/Ismail's argument]\label{ismailarg}
Let  $F(x)$ and $G(x)$ be arbitrary two analytic functions  over the complex field $\mathbb{C}$. If there exists an infinite sequence $\{b_n\}_{n\geq 0}$ with
$\lim_{n\to\infty }b_n=b$, such that for arbitrary $n\geq 0$,
\begin{eqnarray}
F(b_n)=G(b_n),
\end{eqnarray}
then $F(x)=G(x)$ for all $x\in U(b)$, a neighborhood of $b$.
\end{yl}
As briefly mentioned above, the following  $(f,g)$-matrix inversion is central  to our discussion.
\begin{yl}[The $(f,g)$-matrix inversion: {\rm \cite[Theorem 1.3]{0020}}]\label{fglm}
 Let $F=(A_{n,k})_{n,k\in {\mathbb N}}$ and $G=(B_{n,k})_{n,k\in {\mathbb N}}$ be two matrices
with entries given by
\begin{eqnarray}
A_{n,k}&=&\frac{\prod_{i=k}^{n-1}f(x_i,b_k)}
{\prod_{i=k+1}^{n}g(b_i,b_k)}\label{news1115550}\qquad\mbox{and}\\
B_{n,k}&=&
\frac{f(x_k,b_k)}{f(x_n,b_n)}\frac{\prod_{i=k+1}^{n}f(x_i,b_n)}
{\prod_{i=k}^{n-1}g(b_i,b_n)},\label{news1115551}\quad\mbox{respectively},
\end{eqnarray}
where  $\{x_n\}_{n\in {\mathbb N}}$ and
$\{b_n\}_{n\in {\mathbb N}}$ are two arbitrary sequences such that none of the denominators
in the right-hand sides of (\ref{news1115550}) and
(\ref{news1115551}) vanish. Then $F=(A_{n,k})_{n,k\in {\mathbb N}}$ and
$G=(B_{n,k})_{n,k\in {\mathbb N}}$ is a matrix inversion if and
only if for all $a,b,c,x\in \mathbb{{\mathbb C}}$,
 \begin{eqnarray}
g(a,b)f(x,c)+g(b,c)f(x,a)+g(c,a)f(x,b)=0\label{triid}
\end{eqnarray}
with a prior requirement  that  $g(x,y)=-g(y,x).$
\end{yl}

As the earlier work of \cite{0020} displays, there are numerous
functional pairs $f(x,y)$ and $g(x,y)$  satisfying
(\ref{triid}),  each  pair of which in turn constitutes an $(f,g)$-matrix inversion useful in the study of $q$-series.   The reader is referred to \cite{chu, 0020,0021} for further details.

In the present paper, we shall restrict  ourselves to functions  $F(x)$ in (\ref{carbe})-(\ref{chuexp}) which are such a kind of $q$-series that their summands contain a commonly occurring factor $$\frac{\poq{b/x}{n}}{\poq{ax}{n}}\,x^n,$$
where  $a$ and $b$ are two parameters being independent of variable $x$.
Such a factor is now known
 as the  most distinctive feature  of well-poised $q$-series.

With such idea in mind, we shall set up the following    $(1-xy,y-x)$-expansion formula, which is the theoretical   basis for our forthcoming discussions.
\begin{dl}[The $(1-xy,y-x)$-expansion formula] \label{dl4}
Let $F(x)$ be any analytic function over ${\Omega}\subseteq \mathbb{C}$, $\{b_n\}_{n\geq 0}, \{x_n\}_{n\geq 0}\subseteq {\Omega}$ such that
\begin{enumerate}
  \item[(i)] $\{b_n\}_{n\geq 0}$ are pairwise distinct and $\{x_n\}_{n\geq 0}$ is bounded;
  \item[(ii)] $\lim_{n\mapsto \infty}b_n=b\neq b_0$, $\inf\{|1/x_n-b|: n\geq
0\}>0$;
  \item[(iii)] Suppose further that $$\limsup_{n\mapsto
\infty}|G_{n,1}/G_{n,0}|<\infty,$$ where for $\epsilon=0,1$, we define
\begin{align}
G_{n,\epsilon}:=\sum_{k=0}^{n}F(b_{k+\epsilon})\frac{\prod_{i=1}^{n-1} (1-x_{i}b_{k+\epsilon})}
 {\prod_{i=0,i\neq k}^{n}(b_{k+\epsilon}-b_{i+\epsilon})}.\label{coeffepsilon}
\end{align}
\end{enumerate}
 Then  there exists  an open set $\Omega_1\subseteq \Omega$ containing
$b$ such that  for $x\in \Omega_1$,
\begin{align}
F(x)=\sum_{n=0}^{\infty}(1-x_nb_n)G_{n,0}
 \frac{\prod_{i=0}^{n-1}(x-b_i)}{\prod_{i=1}^{n}(1-xx_i)}.\label{1.822}\end{align}
\end{dl}

With the help of  the $(1-xy,y-x)$-expansion formula above and appropriate choices of $F(x)$,   we shall establish

\begin{dl} \label{dlmain-more}With the assumptions given by Theorem \ref{dl4}. Then it holds
\begin{align}
\sum_{n=0}^\infty\alpha_n\frac{\poq{c/x}{n}}{\poq{ax}{n+1}}x^n=\sum_{n=0}^\infty\beta_n\frac{\poq{c/x}{n}}{\poq{bx}{n+1}}x^n,\label{dl4115}
\end{align}
where
\begin{subequations}
\begin{align}
 \alpha_n&=\frac{b^{n}\poq{a/b}{n}(1-acq^{2n})}{\poq{bc}{n+1}}
 \sum_{k=0}^{n}\frac{\poq{bc,q^{-n},acq^{n}}{k}}{\poq{q,bcq^{n+1},bq^{1-n}/a}{k}}(q/a)^{k}\beta_{k},\label{useful-1}\\
 \beta_n&=\frac{a^{n}\poq{b/a}{n}(1-bcq^{2n})}{\poq{ac}{n+1}}\sum_{k=0}^{n}\frac{\poq{ac,q^{-n},bcq^{n}}{k}}{\poq{q,acq^{n+1},aq^{1-n}/b}{k}}(q/b)^{k}\alpha_{k}.\label{useful-2}
\end{align}
\end{subequations}
\end{dl}
From the perspective of applications to basic hypergeometric series,   two special cases of Theorem \ref{dlmain-more} are especially useful.
\begin{tl}\label{dlmain-old} For any integers $s\geq r\geq 0$ and  variables $x,t$, let  $\{a_n\}_{1\leq n\leq r}, \{b_n\}_{1\leq n\leq s},$ and $a,b,c$ be any complex numbers such that all infinite sums converge.  Then it holds
\begin{align}
&{}_{r+1}\widetilde{\phi}_{s+1}\left
[\begin{matrix}a_{1},&a_{2},&\dots,&a_{r},&c/x\\ b_{1},&b_2,&\dots
,&b_s,&bxq\end{matrix} ; q, xt\right]\nonumber\\
&=(1-bx)\sum_{n=0}^{\infty}\frac{(a/b,c/x;q)_{n} }{(bc,ax;q)_{n+1}}(1-acq^{2n})~(bx)^{n}\label{utrans}\\
&\quad\times\,{}_{r+3}\phi _{s+2}\left
[\begin{matrix}bc,&q^{-n},&acq^{n},&a_{1},&\dots,&a_{r}\\ &bcq^{n+1},&bq^{1-n}/a,&b_{1},&\dots
,&b_{s}\end{matrix}; q, \frac{tq}{a}\right].\nonumber
\end{align}
\end{tl}

The limitation of (\ref{utrans}) as $b$ tends to zero together with  $a=1$ and $b_s=q$ deserves our attention, to which we shall frequently have recourse in the sequel
and which we state here as:
\begin{tl}[Expansion of the ${}_{r}\phi_{s-1}$ series in terms of ${}_{r+1}\phi_{s}$ series] \label{dlmain-add}With the same assumption as above. We have
\begin{align}
&{}_{r}\phi _{s-1}\left
[\begin{matrix}a_{1},&a_{2},&\dots,&a_{r-1},&c/x\\ &b_{1},&\dots
,&b_{s-2},&b_{s-1}\end{matrix} ; q, xt\right]\label{dlidi}\\
&=\sum_{n=0}^{\infty}
{}_{r+1}\phi _{s}\left
[\begin{matrix}q^{-n},&a_{1},&\dots,&a_{r-1},&cq^{n}\\ &b_{1},&\dots
,&b_{s-1},&q\end{matrix}; q, tq\right]\frac{(c/x;q)_n }{(x;q)_{n+1}}(1-cq^{2n})q^{n^2/2-n/2}(-x)^n.\nonumber\end{align}
\end{tl}
Furthermore, we can prove
\begin{dl} \label{dlmain}With the assumptions given by Theorem \ref{dl4}.  Then for any integers $s\geq r\geq 0$ and  sequence $\{\beta_n,A_n,B_n\}_{n\geq 0}$ such that all  infinite sums converge, it holds
\begin{align}
&\sum_{n\geq k\geq 0}\frac{1-aq^{2n}}{1-a}\frac{\poq{a,1/x}{n}x^n}{\poq{q,aqx}{n}}\frac{\poq{q^{-n},aq^n,A_1,A_2,\cdots,A_{r}}{k}}{\poq{aq,bq^{k+1},B_1,B_2,\cdots,B_{s}}{k}}
(-1)^kq^{k^2/2+k/2}\beta_k\nonumber\\&\qquad\quad\times{}_{r+3}\phi _{s+2}\left[\begin{matrix}q^{-n+k},&aq^{n+k},bq^{k+1},A_1q^k,A_2q^k,\cdots,A_{r}q^k\\  &aq^{k+1},bq^{2k+1},B_1q^k,B_2q^k,\cdots,B_{s}q^k\end{matrix}
; q, q\right]\label{eq333}\\
&=\frac{\poq{A_1,A_2,\cdots,A_{r},B_1x,B_2x,\cdots,B_{s}x}{\infty}}{\poq{B_1,B_2,\cdots,B_{s},A_1x,A_2x,\cdots,A_{r}x}{\infty}}\sum_{n=0}^\infty \frac{\poq{1/x}{n}}{\poq{bqx}{n}} x^n\beta_n.\nonumber
\end{align}
\end{dl}

All results above shows that one can express any arbitrary modified ${}_{r+1}\phi_{s+1}$ series in variable $x\,t$ as a linear combination of  terminating ${}_{r+3}\phi_{s+2}$ series in variable $t$. As we shall see later, it not only contains the classical Rogers-Fine identity, Carlitz's expansion (\ref{carbe}), Liu's expansion (\ref{1111}), and Chu's expansion (\ref{chuexp}) as special cases in the context of $q$-series, but also draws a general framework for the well-poised Bailey lemma associated with arbitrary well-poised Bailey pair, a concept first appeared  in the work \cite{andrewsber} by Andrews and Berkovich. Besides,   we emphasize here that the $(f,g)$-matrix inversion indeed helps us to overcome the computational difficulty of  the $q$-derivative operator.

 The  remainder of our paper is organized as follows.
In the next section, we shall set up two preliminary results of the $(f,g)$-matrix inversion. With the aid of these results, we shall show Theorems \ref{dl4},\ref{dlmain-more}, \ref{dlmain},  Corollaries \ref{dlmain-old} and \ref{dlmain-add} in full details. Some important specifications of these    conclusions are presented.  In Section 3, we shall apply these  expansion formulas  and allied corollaries  to seek for transformation and summation formulas of $q$-series, among these results are new approaches to the Askey-Wilson polynomials, the Rogers-Fine identity, Andrews'  four-parametric reciprocity theorem and Ramanujan's famous ${}_1\psi_1$ summation formula. Especially noteworthy is  both Theorem \ref{dlmain-more} and Theorem \ref{dlmain}  are closely related to the well-poised Bailey lemma for well-poised Bailey pairs.

 We conclude our introduction with  some remarks on notation.  Throughout this paper, we shall adopt the standard notation and terminology for basic
hypergeometric series (or $q$-series) found in the book  \cite{10} of Gasper and Rahman. Given a
(fixed) complex number $q$ with $|q|<1$, a complex number $a$ and an
integer $n$, define the $q$-shifted factorials
 $(a;q)_\infty$ and $(a; q)_n$  as
\begin{eqnarray*}
(a;q)_{\infty}=\prod_{n=0}^{\infty}(1-aq^n),\quad (a;q)_n=
\frac{\poq{a}{\infty}}{\poq{aq^n}{\infty}}. \label{conven}
    \end{eqnarray*}
We also employ the following compact multi-parameter notation
\begin{eqnarray*}
(a_1,a_2,\cdots,a_m;q)_n = (a_1;q)_n(a_2;q)_n\cdots
   (a_m;q)_n.
    \end{eqnarray*}
The basic and bilateral hypergeometric series with the
base $q$ and  variable $x$   are defined respectively as
\begin{eqnarray*}{}_{r}\phi _{s}\left[\begin{matrix}a_{1},\dots ,a_{r}\\ b_{1},\dots ,b_{s}\end{matrix}
; q, x\right]&=&\sum _{n=0} ^{\infty }\frac{\poq {a_{1},\cdots
,a_{r}}{n}}{\poq
{q,b_{1},\cdots,b_{s}}{n}}\,\tau(n)^{1+s-r}x^{n},\\
{}_{r}\psi _{r}\left[\begin{matrix}a_{1},\dots ,a_{r}\\
b_{1},\dots ,b_{r}\end{matrix}; q, x\right]&=&\sum _{n=-\infty}
^{\infty }\frac{\poq {a_{1},\cdots ,a_{r}}{n}}{\poq
{b_{1},\cdots,b_{r}}{n}}x^{n}.
\end{eqnarray*}
Hereafter, for  convenience we write $\tau(n)$ for $(-1)^nq^{\binom{n}{2}}$. In particular, when $s=r-1$ and the
parameters above satisfy the relations
$$
qa_1=b_1a_2=\cdots=b_{r-1}a_{r},
$$
we  call such a basic hypergeometric
series \emph{well-poised} (in short, WP)  and further, if $a_2=q\sqrt{a_1},a_3=-q\sqrt{a_1},$
 \emph{very-well-poised} (VWP). As is customary, we denote the latter by the compact notation
\[\,_{r}W_{r-1}(a_1;a_4,a_5,\cdots,a_{r};q,x).\]
For the convenience of our discussion, we introduce a modified basic hypergeometric series
$$
{}_{r}\widetilde{\phi} _{s}\left[\begin{matrix}a_{1},\dots ,a_{r}\\ b_{1},\dots ,b_{s}\end{matrix}
; q, x\right]=\sum _{n=0} ^{\infty }\frac{\poq {a_{1},\cdots
,a_{r}}{n}}{\poq
{b_{1},\cdots,b_{s}}{n}}\,\tau(n)^{s-r}x^{n}.
$$
\section{Proofs of the main results and allied corollaries}
\subsection{The $(f,g)$-expansion formula}
Let us begin with  two preliminary results. One is a variation of  the $(f,g)$-matrix inversion, viz., Lemma \ref{fglm}. Another is the $(f,g)$-expansion formula which can be regarded as a direct application of this inversion formula to series expansions for analytic functions.

\begin{yl}\label{dl2}
 Let $\{x_n\}_{n\in \mathbb{N}}$ and $\{b_n\}_{n\in \mathbb{N}}$
 be arbitrary sequences over $\mathbb{C}$ such that $b_n, n\in \mathbb{N}=\{0,1,2,\ldots\}$, are pairwise
 distinct, $g(x,y)=-g(y,x),f(x,y)$ is subject to $(\ref{triid})$.
  Then the linear system with respect to two sequences $\{F_n\}_{n\in \mathbb{N}}$ and $\{G_n\}_{n\in \mathbb{N}}$ \begin{eqnarray}
F_n=\sum_{k=0}^{n}G_k f(x_k,b_k)
\frac{\prod_{i=0}^{k-1}g(b_i,b_n)}{\prod_{i=1}^{k}f(x_i,b_n)}\label{27}
\end{eqnarray}
is equivalent to \begin{eqnarray}
G_n=\sum_{k=0}^{n}F_k\frac{\prod_{i=1}^{n-1} f(x_i,b_k)}
 {\prod_{i=0,i\neq k}^{n}g(b_i,b_k)}.\label{28}
\end{eqnarray}
\end{yl}
\begin{proof}\,\,First,  assume that (\ref{28}) holds for $n\geq 0$.
It can evidently  be manipulated as
$$
\sum_{k=0}^{n}
 \frac{\prod_{i=k}^{n-1}f(x_i,b_k)}
 {\prod_{i=k+1}^{n}g(b_i,b_k)}\left\{\frac{\prod_{i=1}^{k-1}f(x_i,b_k)}
 {\prod_{i=0}^{k-1}g(b_i,b_k)}F_k\right\} =G_n.
$$
 Now, with the  help  of Lemma \ref{fglm}, we can solve this linear system (in infinite unknown members $F_k$) for the terms within the curly brackets.
 The solution is as follows:
$$
\sum_{k=0}^{n}G_k f(x_k,b_k)\frac{\prod_{i=k+1}^{n-1}f(x_i,b_n)}
{\prod_{i=k}^{n-1}g(b_i,b_n)}=
 \frac{\prod_{i=1}^{n-1}f(x_i,b_n)}
 {\prod_{i=0}^{n-1}g(b_i,b_n)}F_n.
$$
 Upon multiplying both sides by $\prod_{i=0}^{n-1}g(b_i,b_n)/\prod_{i=1}^{n-1}f(x_i,b_n)$, we   obtain
$$
F_n=\sum_{k=0}^{n}G_kf(x_k,b_k)
\frac{\prod_{i=0}^{k-1}g(b_i,b_n)}{\prod_{i=1}^{k}f(x_i,b_n)}.
$$
Hence, the identity (\ref{27}) is proved. Conversely, assume that (\ref{27}) is true, from which  one may deduce (\ref{28})  by the same argument. The lemma is thus proved.\end{proof}

\begin{yl}[The $(f,g)$-expansion formula] \label{yl3}  With all conditions of Lemma \ref{dl2}. If there exists an expansion of the form
\begin{eqnarray}
F(x)=\sum_{k=0}^\infty G_k f(x_k,b_k)\frac{\prod_{i=0}^{k-1}g(b_i,x)}{\prod_{i=1}^{k}f(x_i,x)},\label{equfun0}
\end{eqnarray}
then the coefficients
\begin{eqnarray}
G_n=\sum_{k=0}^{n}F(b_k)\frac{\prod_{i=1}^{n-1} f(x_i,b_k)}
 {\prod_{i=0,i\neq k}^{n}g(b_i,b_k)}.\label{280}
\end{eqnarray}
\end{yl}
\begin{proof}
 It only needs to  set $x=b_n$ in (\ref{equfun0}). On account of the fact that $g(b_n,b_n)=0$,
 we immediately come up with a linear system
\begin{eqnarray*}
 F(b_n)=\sum_{k=0}^n G_kf(x_k,b_k)\frac{\prod_{i=0}^{k-1}g(b_i,b_n)}{\prod_{i=1}^{k}f(x_i,b_n)}. \label{equi}
  \end{eqnarray*}
  By virtue of Lemma \ref{dl2}, we thereby obtain (\ref{280}).
\end{proof}

\begin{remark}\label{rmk}
 It is worthwhile mentioning that (\ref{equfun0}) furnishes a series expansion for arbitrary analytic function $F(x)$, provided that the infinite sum on the right-hand side is convergent for  $x, b_n, x_n\in U(b)$ such that $b_n\to b$ and $f(x,x_n)\neq 0$. Some remarkable results of such sort, besides the expansion  (\ref{gessel1})/(\ref{gessel2}) above, are the $q$-Taylor theorems due to Ismail and Stanton \cite{ismail1,ismail2}.
\end{remark}

To make our argument easier to
follow, we record a few basic relations  of $q$-series as follows.

\begin{yl}[\mbox{\cite[cf.][(I.10)-(I.11)]{10}}]\label{basicrelations} Let $\poq{a}{n}$ and $\tau(n)$ be defined as above. Then
\begin{subequations}
\begin{align}
   \frac{\poq{a}{n-k}}{\poq{b}{n-k}}
   &=\frac{\poq{a}{n}}{\poq{b}{n}}\frac{\poq{q^{1-n}/b}{k}}
 {\poq{q^{1-n}/a}{k}}\bigg(\frac{b}{a}\bigg)^k.\label{posineg}\\
 \mbox{In particular,}\quad \poq{a}{n-k}
   &=(-1/a)^kq^{\binom{k+1}{2}-nk}\frac{\poq{a}{n}}
   {\poq{q^{1-n}/a}{k}}\label{posineg1}.\\
   \tau(n+k)&=\tau(n)\tau(k)\,q^{nk}\label{posineg2}.
\end{align}
\end{subequations}
\end{yl}
In what follows,   we shall frequently appeal to the $q$-Pfaff-Saalsch\"{u}tz ${}_3\phi_2$ summation formula.
\begin{yl}[\mbox{\cite[cf.][(II.12)]{10}}]\label{Saalschutz} For integers $n\geq 0$ and $abq^{1-n}=cd$, it holds
\begin{align}
  {}_{3}\phi _{2}\left
[\begin{matrix}q^{-n},&a,&b\\ &c,&d\end{matrix}; q, q\right]=\frac{\poq{c/a,c/b}{n}}{\poq{c,c/(ab)}{n}}\label{saal}.
\end{align}
\end{yl}

%\subsection{Proofs of the main results}
\subsection{The analytic proof of the $(1-xy,y-x)$-expansion formula}
Now we are ready to show Theorem \ref{dl4}.

\begin{proof}\,\, Let $r_0=\inf\{|1/x_n-b|: n\geq 0\}$. The
assumption that $r_0>0$ means that  each function $1/(1-xx_i)$ has
no pole in the disk ${\Bbb O}_{r_0}=\{x:|x-b|<r_0\}$.  For
notational simplicity, we write
\begin{eqnarray}S_n(x)=\sum_{k=0}^{n}G_{k,0}(1-x_kb_k)
\frac{\prod_{i=0}^{k-1}(x-b_i)}{\prod_{i=1}^{k}(1-xx_i)}.\label{expandef1}
\end{eqnarray}
It is of importance to see that $S_{n}(b_n)=F(b_n)$.
Thus, according
to Lemma \ref{ismailarg}, it needs only to show that there exists
an open set $\Omega_1$ containing $b$, such that $S_n(x)\mapsto
S_{\infty}(x)$ uniformly on $\Omega_1\subseteq {\Bbb O}_{r_0}$.   To this end, it is easy to check that for the coefficients $G_{n,\epsilon}$ given by \eqref{coeffepsilon}, it holds
 the recursive formula as follows:
\begin{align}
G_{k+1,0}=
\frac{1-x_kb_0}{b_0-b_{k+1}}G_{k,0}+\frac{1-x_kb_{k+1}}{b_{k+1}-b_0}G_{k,1}.
\end{align}
Evidently, the ratio
 \begin{align*}
\left|\frac{G_{k+1,0}}{G_{k,0}}\right|&=\left|\frac{1-x_kb_0}{b_0-b_{k+1}}+\frac{G_{k,1}}{G_{k,0}}\frac{1-x_kb_{k+1}}{b_{k+1}-b_0}\right|\\
&=\bigg|\frac{1-x_kb_0}{b_0-b_{k+1}}+\lambda_k\frac{1-x_kb_{k+1}}{b_{k+1}-b_0}\bigg|
=\bigg|-x_k+(\lambda_k-1)\frac{1-x_kb_{k+1}}{b_{k+1}-b_0}\bigg|.\end{align*}
From here, by the triangular inequality, we have
 \begin{align}
\left|\frac{G_{k+1,0}}{G_{k,0}}\right|\leq
|x_k|+|\lambda_k-1|\left|\frac{1-x_kb_{k+1}}{b_{k+1}-b_0}\right|.\label{budengshi-old}
\end{align}
 Hence, because
$\lim_{k\mapsto\infty}b_k=b$, we see that for given $\epsilon,$ there exists an
integer $K_1$ that for $k>K_1$, such that $|b_k-b|<\epsilon$. Therefore, we have
\begin{align}
\left|\frac{1-x_kb_{k+1}}{b_{k+1}-b_0}\right|=\left|\frac{1-x_kb_{0}}{b_{k+1}-b_0}-x_k\right|\leq\frac{1+|x_kb_{0}|}{|b_{k+1}-b_0|}+|x_k|<M_1:= m+\frac{1+|b_0|m}{|b-b_0|-\epsilon},\label{budengshi}
\end{align}
 where $m=\sup\{|x_k|:k\geq 0\}$.  With \eqref{budengshi-old} and  \eqref{budengshi} taken into account, it is easy to see that
for $k>K>K_1$,
  \begin{eqnarray}
\left|\frac{G_{k+1,0}}{G_{k,0}}\right|\leq
|x_k|+|\lambda_k-1|\left|\frac{1-x_kb_{k+1}}{b_{k+1}-b_0}\right|<m+m_0M_1=M,\label{inequ1}
\end{eqnarray}
where $m_0=1+\limsup_{k\mapsto \infty}|\lambda_k|$.
It remains to decide a constant $s<r_0$  by solving the following
inequality
\begin{align}
\left|\frac{1-x_{k+1}b_{k+1}}{1-x_{k}b_{k}}\frac{b_k-x}{1-x_{k+1}x}\right|<\frac{a}{M}\left(<\frac{1}{m}\right)\label{inequ2}\end{align}
for the fixed $M$ and $a:0<a<1$,  such that $|x-b|<s$ small enough while $s$ is
independent of $k$ and $x$.
Consequently, the inequalities (\ref{inequ1}) and (\ref{inequ2}) together asserts that for $k>K$,
\begin{eqnarray*}
\left|\frac{\mbox{$(k+1)$-th term  of (\ref{1.822})}}{\mbox{$k$-th term of
(\ref{1.822})}}\right|=\left|\frac{G_{k+1,0}}{G_{k,0}}\frac{1-x_{k+1}b_{k+1}}{1-x_{k}b_{k}}\frac{b_k-x}{1-x_{k+1}x}\right|
<M\times \frac{a}{M}=a<1.\end{eqnarray*}
 By the Weierstrass M-test \cite[cf.][p.37]{ahl}, we conclude that
$S_n(x)$ indeed converges uniformly on the disk ${\Bbb
O}_s=\{x:|x-b|<s\}$, which is such a required open set $\Omega_1$.
Hence, by Ismail's argument or Lemma \ref{ismailarg}, we see the limit function $S_{\infty}(x)=F(x)$ for $x\in\Omega_1$.  Theorem \ref{dl4} is proven. \end{proof}

\subsection{Proof and WP Bailey pairs of Theorem \ref{dlmain-more}}
Now we  proceed to the proof of Theorem \ref{dlmain-more}.
\begin{proof}It suffices to consider
$$
F(x)=(1-ax)\sum_{n=0}^\infty\beta_n\frac{\poq{c/x}{n}}{\poq{bx}{n+1}}x^n.
$$
Then expand $F(x)$  in terms of $\{\poq{c/x}{n}x^n/\poq{aqx}{n}\}_{n\geq 0}$  with respect to the parametric sequences
\begin{eqnarray}
x_n=aq^{n}, b_n=cq^{n}\,\, \big(\lim_{n\to\infty}b_n=0\big).\label{parameters-new}
 \end{eqnarray}
 Under this situation, by employing Theorem \ref{dl4}  and Lemma \ref{basicrelations}, we can  easily check that $\alpha_0=F(c)$, and for $n\geq 1$,
\begin{align*}
 \frac{\alpha_n}{1-acq^{2n}}
&=\frac{(ac;q)_{n}}{c^n\poq{q}{n}}\sum_{i=0}^n
\frac{\poq{q^{-n},acq^n}{i}}{\poq{ac}{i}\poq{bcq^i}{i+1}}\beta_i(cq)^i\tau(i)T_{n,i},
\end{align*}
where  the inner sum, by virtue of  the $q$-Pfaff-Saalsch\"{u}tz ${}_3\phi_2$  summation formula \eqref{saal}, is  easily found to be
$$
T_{n,i}:=\sum_{k=0}^{n-i}
\frac{\poq{q^{-n+i},acq^{n+i},bcq^i}{k}}{\poq{q,acq^{i},bcq^{2i+1}}{k}}q^k=\frac{\poq{q^{-n},a/b}{n-i}}{\poq{acq^{i},1/(bcq^{n+i})}{n-i}},
$$
giving rise to
\begin{align*}
 \frac{\alpha_n}{1-acq^{2n}}
&=\frac{1}{c^n\poq{q}{n}}\sum_{i=0}^n
\frac{\poq{q^{-n},acq^n}{i}}{\poq{bcq^i}{i+1}}\frac{\poq{q^{-n},a/b}{n-i}}{\poq{1/(bcq^{n+i})}{n-i}}\beta_i(cq)^i\tau(i)\\
&=\frac{b^{n}\poq{a/b}{n}}{\poq{bc}{n+1}}\sum_{i=0}^{n}\frac{\poq{bc,q^{-n},acq^{n}}{i}}{\poq{q,bcq^{n+1},bq^{1-n}/a}{i}}(q/a)^{i}\beta_{i}.
\end{align*}
Hence, \eqref{useful-1} is confirmed. Note that \eqref{useful-2} is obtainable by exchanging  $a$ and $b$. This completes the proof of the theorem.
\end{proof}
Theorem \ref{dlmain-more} brings us to Andrews' the concept of   WP Bailey pairs, which is a generalization of the usual  Bailey pairs.
\begin{dy}[{\rm \cite[cf.][Definition 6.1]{andrews5}}]
  A pair of sequences $\{\alpha_n(t,b)\}_{n\geq 0}$
and $\{\beta_n(t,b)\}_{n\geq 0}$ satisfying
\begin{align}
\,\beta_n(t,b)=\sum^{n}_{k=0} \frac{\poq{b}{n-k}\poq{bt}{n+k}}
{\poq{q}{n-k}\poq{tq}{n+k}}\,\alpha_k(t,b)\label{def}
\end{align}
is called to be a well-poised (WP) Bailey pair with respect to the parameters $t$ and $b$.
\end{dy}
From this viewpoint,  Theorem \ref{dlmain-more} can be recognized as not only a special form of the WP Bailey lemma \cite[Theorem 3.2]{liuqq} but also a  machine to produce  WP Bailey pairs. We state such machine  by the following
\begin{tl}Let  $\{\alpha_n\}_{n\geq 0}$
and $\{\beta_n\}_{n\geq 0}$
 be a  pair of sequences subject to \eqref{useful-1}/\eqref{useful-2}. Then
  \begin{align}\left\{
  \begin{array}{ll}\displaystyle
 \alpha_n(bc,a/b)&:=\displaystyle\frac{\poq{bc}{n}}{\poq{q}{n}}(1/b)^{n}\beta_{n}\\
 &\\
\beta_n(bc,a/b)&:=\displaystyle \frac{\poq{ac}{n}(1-bc)}{\poq{q}{n}(1-acq^{2n})}(1/b)^n\alpha_n
\end{array}
\right.
 \end{align}
 is a  WP Bailey pair.
\end{tl}
\begin{proof} It suffices to restate \eqref{useful-1} in the form
\begin{align*}
 \frac{\poq{bc}{n+1}\alpha_n}{b^{n}\poq{a/b}{n}(1-acq^{2n})}=\sum_{k=0}^{n}\frac{\poq{q^{-n},acq^{n}}{k}}{\poq{bq^{1-n}/a,bcq^{n+1}}{k}}\frac{\poq{bc}{k}}{\poq{q}{k}}(q/a)^{k}\beta_{k},
\end{align*}
or equivalently,
\begin{align*}
  \frac{\poq{ac}{n}\poq{bc}{n+1}\alpha_n}{b^{n}\poq{bcq}{n}\poq{a/b}{n}(1-acq^{2n})}=\sum_{k=0}^{n}\frac{\poq{q^{-n}}{k}\poq{ac}{n+k}}{\poq{bq^{1-n}/a}{k}\poq{bcq}{n+k}}\frac{\poq{bc}{k}}{\poq{q}{k}}(q/a)^{k}\beta_{k}.
\end{align*}
Recall that
$$\frac{\poq{q^{-n}}{k}}{\poq{bq^{1-n}/a}{k}}
=\frac{\poq{a/b}{n-k}}{\poq{q}{n-k}}\frac{\poq{q}{n}}{\poq{a/b}{n}}(a/(bq))^k,
$$
which allows us to simplify further
\begin{align*}
   \frac{\poq{ac}{n}(1-bc)}{\poq{q}{n}(1-acq^{2n})}(1/b)^n\alpha_n&=\sum_{k=0}^{n}\frac{\poq{a/b}{n-k}\poq{ac}{n+k}}{\poq{q}{n-k}\poq{bcq}{n+k}}\frac{\poq{bc}{k}}{\poq{q}{k}}(1/b)^{k}\beta_{k}.
\end{align*}
In view of the definition of  WP Bailey pairs,  the desired conclusion follows.
\end{proof}

Accordingly, the WP Bailey lemma   \cite[cf.][Theorem 3.2]{liuqq} reduces to

\begin{tl}\label{wpbaileylemma} Let  $\{\alpha_n\}_{n\geq 0}$
and $\{\beta_n\}_{n\geq 0}$
 be a  pair of sequences subject to \eqref{useful-1}/\eqref{useful-2}. Then it holds
\begin{align}
&\sum_{n=0}^{\infty}\frac{\poq{ac,x,y}{n}}
{\poq{q,acq/x,acq/y}{n}}\left(\frac{cq}{xy}\right)^n\alpha_n\label{wpbaileylemma-id}\\
&= \frac{\poq{ac,bcq/x,bcq/y,acq/(xy)}{\infty}}{\poq{bc,bcq/(xy),acq/x,acq/y}{\infty}}\sum_{n=0}^{\infty}\frac{\poq{bc,x,y}{n}}
{\poq{q,bcq/x,bcq/y}{n}}\left(\frac{cq}{xy}\right)^n\beta_n
.\nonumber
 \end{align}
\end{tl}

\subsection{Proofs of Corollaries \ref{dlmain-old} and \ref{dlmain-add} }

\begin{proof}[\sl The proof of Corollary \ref{dlmain-old}.]
It suffices to take in \eqref{dl4115} of Theorem \ref{dlmain-more}
\begin{align}
\beta_n:=\tau(n)^{s-r}t^n \frac{\poq{a_1,a_2,\cdots,a_{r}}{n}}{\poq{b_1,b_2,\cdots,b_{s}}{n}},
\end{align}
 so that we can expand
\[F(x)=\frac{1}{1-bx}\sum_{n=0}^\infty \frac{\poq{a_1,a_2,\cdots,a_{r},c/x}{n}}{\poq{b_1,b_2,\cdots,b_{s},bqx}{n}}\tau(n)^{s-r}(xt)^n \]
 into
 \begin{align*}
\sum_{n=0}^\infty\alpha_n\frac{\poq{c/x}{n}}{\poq{ax}{n+1}}x^n.
\end{align*}
By \eqref{useful-1} of Theorem \ref{dlmain-more}, we obtain
\begin{align*}
\alpha_n&=\frac{b^{n}\poq{a/b}{n}(1-acq^{2n})}{\poq{bc}{n+1}}\sum_{i=0}^{n}\frac{\poq{bc,q^{-n},acq^{n},a_1,a_2,\cdots,a_{r}}{i}}{\poq{q,bcq^{n+1},bq^{1-n}/a,b_1,b_2,\cdots,b_{s}}{i}}\tau(i)^{s-r}(tq/a)^{i}\\
&=\frac{b^{n}\poq{a/b}{n}(1-acq^{2n})}{\poq{bc}{n+1}}\,{}_{r+3}\phi _{s+2}\left
[\begin{matrix}bc,&q^{-n},&acq^{n},&a_{1},&\dots,&a_{r}\\ &bcq^{n+1},&bq^{1-n}/a,&b_{1},&\dots
,&b_{s}\end{matrix}; q, \frac{tq}{a}\right].
\end{align*}
In the end we achieve
\begin{align*}
&{}_{r+1}\widetilde{\phi}_{s+1}\left
[\begin{matrix}a_{1},&a_{2},&\dots,&a_{r},&c/x\\ b_{1},&b_2,&\dots
,&b_s,&bxq\end{matrix} ; q, xt\right]
=(1-bx)\sum_{n=0}^{\infty}(1-acq^{2n})~(bx)^{n}\frac{(a/b,c/x;q)_{n} }{(bc,ax;q)_{n+1}}\\
&\qquad\qquad\qquad\qquad\qquad\times\,{}_{r+3}\phi _{s+2}\left
[\begin{matrix}bc,&q^{-n},&acq^{n},&a_{1},&\dots,&a_{r}\\ &bcq^{n+1},&bq^{1-n}/a,&b_{1},&\dots
,&b_{s}\end{matrix}; q, \frac{tq}{a}\right].\nonumber
\end{align*}
That is we wanted.
\end{proof}

\begin{proof}[\sl The proof of Corollary \ref{dlmain-add}.]
To establish this expansion, we only need to take the limit  of  both sides of (\ref{utrans}) as $b$ tends to zero and then let $a=1,a_{r}=0,$ and $b_{s}=q$.
\end{proof}

 In what follows, we shall consider a few specific cases of Corollaries \ref{dlmain-old} and \ref{dlmain-add} which will be frequently used. At first,  when  $c$ tends to zero, the expansion (\ref{dlidi}) in Corollary \ref{dlmain-add} furnishes a generalization of Carlitz's $q$-expansion formula for $q$-series \cite[cf.][]{carlitz}.

\begin{tl}\label{tl333} With the same assumptions of  Corollary \ref{dlmain-old}. Then
\begin{eqnarray}
&&{}_{r}\phi _{s}\left
[\begin{matrix}a_{1},&a_2,&\dots,&a_{r-1},&0\\ &b_{1},&\dots
,&b_{s-1},&b_s\end{matrix} ; q, xt\right]\label{utrans-II}\\
&=&\sum_{n=0}^{\infty}{}_{r+1}\phi _{s+1}\left
[\begin{matrix}q^{-n},&a_{1},&\dots,&a_{r-1},&0\\ &b_{1},&\dots
,&b_{s},&q\end{matrix}; q, tq\right]\frac{q^{n^2/2-n/2}(-x)^n}{(x;q)_{n+1}}\,\nonumber
.\end{eqnarray}
\end{tl}

Once putting $r=s$ and $a=1$ in Corollary \ref{dlmain-old}, we obtain a previously unknown transformation  which deserves our particular attention.

\begin{tl}\label{tl4} With the same assumption of  Corollary \ref{dlmain-old}. Then
\begin{align}
&{}_{r+1}\widetilde{\phi}_{r+1}\left
[\begin{matrix}a_{1},&a_2,&\dots,&a_{r},&c/x\\ b_{1},&b_2,&\dots,&b_{r},&bqx\end{matrix}; q, xt\right]\label{utrans-I}\\
&=(1-bx)\sum_{n=0}^{\infty}\frac{(1/b,c/x;q)_{n} }{(bc,x;q)_{n+1}}~(1-cq^{2n})(bx)^{n}
{}~_{r+3}\phi _{r+2}\left
[\begin{matrix}bc,&q^{-n},&cq^{n},&a_{1},&\dots,&a_{r}\\ &bcq^{n+1},&bq^{1-n},&b_{1},&\dots,&b_{r}\end{matrix}; q, tq\right].\nonumber\end{align}
\end{tl}

Assuming further $b_r=q$, then Corollary \ref{tl4}  can be reformulated in the form
\begin{tl}[Expansion of the ${}_{r+1}\phi_{r}$ in terms of the ${}_{r+3}\phi_{r+2}$ series] \label{tl5} With the same assumption of  Corollary \ref{dlmain-old}. Then
\begin{align}
{}_{r+1}\phi _{r}\left
[\begin{matrix}a_1,&a_2,&\dots,&a_{r},&c/x\\ &b_{1},&\dots,&b_{r-1},&bqx\end{matrix}; q, xt\right]\nonumber\\
=(1-bx)\sum_{n=0}^{\infty}
\frac{(1/b,c/x;q)_{n} }{(bc,x;q)_{n+1}}~(1-cq^{2n})(bx)^{n}\\
\times{}~_{r+3}\phi _{r+2}\left
[\begin{matrix}bc,&q^{-n},&cq^{n},&a_{1},&\dots,&a_{r-1},&a_{r}\\ &bcq^{n+1},&bq^{1-n},&b_{1},&\dots,&b_{r-1},&q\end{matrix}; q, tq\right].\nonumber\end{align}
\end{tl}
Of all instances covered by Corollary \ref{tl5}, the most interesting one is that $$a_r=bc; a_1=q\sqrt{bc},a_2=-q\sqrt{bc};\,\,
 a_ib_i=bcq\quad (1\leq i\leq r-1),$$ leading to the following transformation of two VWP series.

 \begin{tl}[Expansion of the ${}_{r+1}W_{r}$  series in terms of the ${}_{r+3}W_{r+2}$ series] \label{tl8}
  Under the same conditions as above. Then
\begin{align}
&{}_{r+1}W_{r}(bc;a_3,a_4,\cdots,a_{r-1},c/x;q,xt)\label{fffiiinnn}\\ &=(1-bx)\sum_{n=0}^{\infty}
\frac{(1/b,c/x;q)_{n} }{(bc,x;q)_{n+1}}~(1-cq^{2n})(bx)^{n}{}~_{r+3}W_{r+2}(bc;q^{-n},cq^{n+1},a_{3},a_4,\dots,a_{r-1},bc; q, tq).\nonumber\end{align}
\end{tl}
\subsection{Proof  and transformations of Theorem \ref{dlmain}}
In this part we begin  to show Theorem \ref{dlmain}, from which we then derive a few useful transformations.
\begin{proof} For convenience, let denote the right-hand side
of \eqref{eq333} by
\begin{align}
F(x)=K(x)\sum_{n=0}^\infty \beta_n \frac{\poq{1/x}{n}}{\poq{bqx}{n}}x^n,
\end{align}
where
\[K(x)=\frac{\poq{A_1,A_2,\cdots,A_{r}}{\infty}}{\poq{B_1,B_2,\cdots,B_{s}}{\infty}}\frac{\poq{B_1x,B_2x,\cdots,B_{s}x}{\infty}}{\poq{A_1x,A_2x,\cdots,A_{r}x}{\infty}}.\]
Clearly,
\[K(q^k)=\frac{\poq{A_1,A_2,\cdots,A_{r}}{k}}{\poq{B_1,B_2,\cdots,B_{s}}{k}}.\]
Suppose that
$$
F(x)=\sum_{n=0}^\infty \alpha_n \frac{\poq{1/x}{n}}{\poq{aqx}{n}}x^n.
$$Thus, according to  Theorem \ref{dl4}, we have
\begin{align*}
 \alpha_n
&=\frac{1-aq^{2n}}{1-a}\frac{(a;q)_{n}}{\poq{q}{n}}\sum_{k=0}^nq^k
\frac{\poq{q^{-n},aq^n}{k}}{\poq{q,aq}{k}}K(q^k)\sum_{i=0}^k\beta_i\frac{\poq{q^{-k}}{i}}{\poq{bq^{k+1}}{i}}q^{ki}\\
&=\frac{1-aq^{2n}}{1-a}\frac{(a;q)_{n}}{\poq{q}{n}}
\sum_{i=0}^n\frac{\poq{q^{-n},aq^n,A_1,A_2,\cdots,A_{r}}{i}}{\poq{aq,bq^{i+1},B_1,B_2,\cdots,B_{s}}{i}}\beta_iq^i\tau(i)U_{n,i},
\end{align*}
where
\begin{align*}
U_{n,i}&:=\sum_{j=k-i=0}^{n-i}
\frac{\poq{q^{-n+i},aq^{n+i},bq^{i+1},A_1q^i,A_2q^i,\cdots,A_{r}q^i}{j}}{\poq{q,aq^{i+1},bq^{2i+1},B_1q^i,B_2q^i,\cdots,B_{s}q^i}{j}}q^j\\
&={}_{r+3}\phi _{s+2}\left[\begin{matrix}q^{-n+i},&aq^{n+i},bq^{i+1},A_1q^i,A_2q^i,\cdots,A_{r}q^i\\  &aq^{i+1},bq^{2i+1},B_1q^i,B_2q^i,\cdots,B_{s}q^i\end{matrix}
; q, q\right].
\end{align*}
So \eqref{eq333} follows.
\end{proof}
From Theorem \ref{dlmain} it follows that
\begin{tl}\label{tl213} For any integers $s\geq r\geq 0$ and any sequence $\{A_n\}_{n\geq 0}$ and $\{B_n\}_{n\geq 0}$ such that all  infinite sums converge, it holds
\begin{align}
&\sum_{n=0}^\infty{}_{r+2}\phi _{s+1}\left[\begin{matrix}q^{-n},&aq^{n},A_1,A_2,\cdots,A_{r}\\  &aq,B_1,B_2,\cdots,B_{s}\end{matrix}
; q, q\right]\frac{1-aq^{2n}}{1-a}\frac{\poq{a,1/x}{n}}{\poq{q,aqx}{n}}x^n\nonumber\\
&\qquad=\frac{\poq{A_1,A_2,\cdots,A_{r},B_1x,B_2x,\cdots,B_{s}x}{\infty}}{\poq{B_1,B_2,\cdots,B_{s},A_1x,A_2x,
\cdots,A_{r}x}{\infty}}.\label{eq334}
\end{align}
\end{tl}
\begin{proof} Observe that \eqref{eq334} is the direct consequence  of \eqref{eq333} when $\beta_n=\delta_{n,0}$.
\end{proof}

Especially noteworthy is that the special case $r=s=3$ of \eqref{eq334} leads us to Liu's extension of Rogers' VWP ${}_6\phi_5$ summation formula \cite{liunew}.  Liu obtained this result via the $q$-exponential differential operator. The reader may consult \cite{liunew2} for details.
\begin{llz}[{\rm \cite[cf.][Theorem 1.5]{liunew}}]  For $|abcd/q^2|<1$, it holds
\begin{align}
&\sum_{n=0}^\infty\frac{1-aq^{2n}}{1-a}\frac{\poq{a,q/b,q/c,q/d}{n}}{\poq{q,ab,ac,ad}{n}}{}_{4}\phi _{3}\left
[\begin{matrix}q^{-n},aq^n,\beta,\gamma\\ q/b,q/c,\beta\gamma\,abc/q\end{matrix}; q, q\right]\bigg(\frac{abcd}{q^2}\bigg)^n\nonumber\\
&=\frac{\poq{aq,\beta\,abc/q,\gamma\,abc/q,abd/q,acd/q,\beta\gamma\,abcd/q^2}{\infty}}
{\poq{ab,ac,ad,\beta\gamma\,abc/q, \beta\,abcd/q^2,\gamma\,abcd/q^2}{\infty}}.
\end{align}
\end{llz}
\begin{proof} It suffices to set
 $r=s=3$ and make the specialization of the parameters
\bnm
(A_1,A_2,A_3)&\to& (aq,\beta\,abc/q,\gamma\,abc/q),\\
x&\to& d/q,\\
(B_1,B_2,B_3)&\to& (ab,ac,\beta\gamma\,abc/q)
\enm
 in \eqref{eq334}.  As such,  it is clear that the ${}_5\phi_4$ series on the left-hand side
\begin{align*}
&{}_{5}\phi _{4}\left
[\begin{matrix}q^{-n},aq^n,aq,\beta\,abc/q,\gamma\,abc/q\\ aq,ab,ac,\beta\gamma\,abc/q\end{matrix}; q, q\right]\\
&={}_{4}\phi _{3}\left
[\begin{matrix}q^{-n},aq^n,\beta\,abc/q,\gamma\,abc/q\\ ab,ac,\beta\gamma\,abc/q\end{matrix}; q, q\right]=\frac{\poq{q/b,q/c}{n}}{\poq{ab,ac}{n}}{}_{4}\phi _{3}\left
[\begin{matrix}q^{-n},aq^n,\beta,\gamma\\ q/b,q/c,\beta\gamma\,abc/q\end{matrix}; q, q\right]\bigg(\frac{abc}{q}\bigg)^n.
\end{align*}Note that the last equality is based on  Sears' transformation \cite[cf.][(III.15)]{10}. This reduces \eqref{eq334} to the claimed.
\end{proof}
There is one general transformation deserving our attention.  As an illustration of  Theorem  \ref{dlmain},  it is not only a generalization  of Rogers'  VWP ${}_6\phi_5$ summation formula \cite[][(II.20)]{10},  it is  also equivalent to the WP Bailey lemma \eqref{wpbaileylemma-id}.
\begin{tl}\label{liu-trnasf} For  any sequence $\{\beta_n\}_{n\geq 0}$ such that all  infinite sums converge, it holds
\begin{align}
&\sum_{n=0}^\infty\frac{1-aq^{2n}}{1-a}\frac{\poq{a,1/x,a/b,aq/B_2}{n}}{\poq{q,aqx,bq,B_2}{n}}(A_2x)^n\alpha_n
\label{good}\\
&=\frac{\poq{aq,bqx,A_2,B_2x}{\infty}}{\poq{bq,aqx,B_2,A_2x}{\infty}}\sum_{n=0}^\infty \frac{\poq{1/x}{n}}{\poq{bqx}{n}} x^n\beta_n,\nonumber
\end{align}
where $aA_2=bB_2$ and
\begin{align}
\alpha_n:=\sum_{k=0}^n\frac{\poq{q^{-n},aq^n,A_2}{k}}{\poq{bq^{n+1},q^{1-n}b/a,aq/B_2}{k}}(q/B_2)^k\beta_k.\label{jacksonid}
\end{align}
\end{tl}
\begin{proof}
It suffices to specialize Theorem \ref{dlmain} to $r=s=2$ and $(A_1,B_1)\to(aq,bq),\quad aA_2=bB_2.
$
In the sequel, we readily find that the ${}_5\phi_4$ series on the left-hand side of \eqref{eq333}, namely,
\begin{align*}
&{}_{5}\phi _{4}\left[\begin{matrix}q^{-n+k},aq^{n+k},bq^{k+1},aq^{k+1},A_2q^k\\  aq^{k+1},bq^{2k+1},
bq^{k+1},B_2q^k\end{matrix}
; q, q\right]\\
&={}_{3}\phi _{2}\left[\begin{matrix}q^{-n+k},aq^{n+k},A_2q^k\\  bq^{2k+1},B_2q^k\end{matrix}
; q, q\right]=\frac{\poq{q^{k+1-n}b/a,q^{k+1}b/A_2}{n-k}}{\poq{bq^{2k+1},q^{1-n}/B_2}{n-k}}.
\end{align*}
The last equality results from the $q$-Pfaff-Saalsch\"{u}tz ${}_3\phi_2$ summation formula \eqref{saal}. This reduces \eqref{eq333} to the desired identity.
\end{proof}
 Indeed, when $\beta_n=\delta_{n,0}$, then \eqref{good} reduces to  Rogers'  VWP ${}_6\phi_5$ summation formula \cite[][(II.20)]{10}, while once specialized to $(x,B_2)\to (1/x,aq/y)$, it turns out to be  the WP Bailey lemma \eqref{wpbaileylemma-id}. In this sense, Corollary \ref{liu-trnasf} may be regarded as a variant of the WP Bailey lemma. Next is a more specific example of \eqref{good}.
 \begin{llz} Let  $aA_2=bB_2$. Then it holds
\begin{align}
\sum_{n=0}^\infty\frac{1-aq^{2n}}{1-a}\frac{\poq{a,c,a/b,aq/B_2}{n}}{\poq{q,aq/c,bq,B_2}{n}}&{}_{4}\phi _{3}\left[\begin{matrix}q^{-n},aq^n,A_2,d\\ bq^{n+1},q^{1-n}b/a,aq/B_2,\end{matrix}
; q, \frac{bq^2}{cdB_2}\right](A_2/c)^n\nonumber\\
&=\frac{\poq{aq,A_2,B_2/c,bq/c^2,bq/(cd)}{\infty}}{\poq{bq,B_2,aq/c,bq/(c^2d),A_2/c}{\infty}}.
\end{align}
\end{llz}
\begin{proof}
It suffices to replace $x$ with $1/c$ and specialize  \eqref{good} to $$ \beta_n=\frac{\poq{d}{n}}{\poq{q}{n}}\bigg(\frac{bq}{cd}\bigg)^n.
$$
By Gauss' ${}_2\phi_1$ summation formula, we readily find that  the series  on the right-hand side of \eqref{good} equals
\begin{align*}
{}_{2}\phi _{1}\left[\begin{matrix}c,&d\\  &bq/c\end{matrix}
; q, \frac{bq}{c^2d}\right]=\frac{\poq{bq/c^2,bq/(cd)}{\infty}}{\poq{bq/c,bq/(c^2d)}{\infty}}.
\end{align*}
It gives the desired identity.
\end{proof}
To our surprise, Corollary \ref{liu-trnasf} contains a transformation of \cite{liunew} by Liu.  Liu used his result to show an important $q$-beta integral with twelve parameters.  Very recently,   Chen and Wang \cite{chenwang} applied Liu's transformation to study representations of mock theta functions.  Our argument for this may go as follows.
\begin{llz}[{\rm \cite[cf.][Proposition 4.1]{liunew}}]  We have
\begin{align}
\sum_{n=0}^\infty\frac{1-aq^{2n}}{1-a}\frac{\poq{a,q/A,q/B}{n}}{\poq{q,aA,aB}{n}}&(-aAB/q)^nq^{n(n-1)/2}
\sum_{k=0}^n\frac{\poq{q^{-n},aq^n}{k}}{\poq{q/A}{k}}(q^2/A)^k\beta_k\nonumber
\label{good-1}\\
&=\frac{\poq{aq,aAB/q}{\infty}}{\poq{aA,aB}{\infty}}\sum_{n=0}^\infty \poq{q/B}{n} (aB)^n\beta_n.
\end{align}
\end{llz}
\begin{proof}
It suffices to specialize   \eqref{good} to the case
$$(A_2,\beta_k)\to(bA,(aq)^k\beta_k).$$ Then let $b$ tend to zero in the identity resulted  by noting that $B_2=aA$. After the replacement of $x$ by $B/q$, we obtain \eqref{good-1}.
\end{proof}
\section{Applications}
\setcounter{equation}{0}
This whole section is devoted to applications of all theorems and corollaries obtained in the preceding section in finding of summation and transformation formulas for $q$-series.
\subsection{A new proof of Rogers' VWP ${}_6\phi_5$ summation formula}
For our purpose, let assume that all parameters in (\ref{utrans}) are independent of $t$. Then both sides of (\ref{utrans}) are power series in $t$. This, together with Ismail's argument, gives  us a new way to Rogers's VWP ${}_6\phi_5$ summation formula \cite[cf.][(II.20)]{10}.
\begin{lz}\label{dladded} Let $a,b,c,d,x: |bx|<1 $ be any complex numbers such that none of the denominators
on both sides of (\ref{tl66655}) below vanish.  Then
\begin{eqnarray}
{}_{6}W_{5}(acd^2;d,aq/b,cd/x;q, bx)=\frac{\poq{axq,bcd,bdx,acd^2q}{\infty}}
{\poq{adxq,bcd^2,bx,acdq}{\infty}}.\label{tl66655}
\end{eqnarray}
\end{lz}
\begin{proof} Under the assumption that each of $a_i,b_j~(1\leq i\leq r,1\leq j\leq s),$ and $x,b,c$ is independent of $t$, we can compare the coefficients of $t^m$  on both sides of (\ref{utrans}) and obtain
\begin{eqnarray}
{}_{6}W_{5}(acq^{2m+2};q^{m+1},aq/b,cq^{m+1}/x; q, bx)=\frac{\poq{axq,bcq^{m+1}}{m+1}}{\poq{acq^{m+2},bx}{m+1}}.\label{tl66655-1}
\end{eqnarray}
Let us denote the sum on the left-hand side of (\ref{tl66655}) by $F(d)$ and the product on the right-hand side by $G(d)$, respectively. Then it is easily found that (\ref{tl66655-1}) amounts to
\[F(q^{m+1})=G(q^{m+1})\quad\mbox{for}\,\,  m\geq 0.\]
 At this stage, it is important to note that both $F(d)$ and $G(d)$ are analytic functions of  variable $d$. Hence, according to Ismail's argument, we conclude that
\[F(d)=G(d)\quad (d\in U(0)).\]
 By analytic continuation, this equality  remains valid for all $a,b,c,d,x$ with $|bx|<1$, provided that the denominators
in both sides of (\ref{tl66655}) are not zero.
\end{proof}
\subsection{Generating functions for the Askey-Wilson polynomials}
From Corollary \ref{dlmain-add} it is easy to derive  a new generating function for the Askey-Wilson polynomials $p_n(y;a,b,c,d|q)$  \cite[cf.][]{10,mar}.
\begin{lz}Let $y=\cos(\theta)$. Then, for $|x|<1$, it holds
\begin{eqnarray}
&&(1-x){}_{3}\widetilde{\phi} _{3}\left
[\begin{matrix}ae^{i\theta},&ae^{-i\theta},&abcd/(xq)\\ ab,&ac
,&ad\end{matrix} ; q, x\right]\label{wilson}\\
&=&\sum_{n=0}^{\infty}\frac{(abcd/(xq);q)_n }{(xq,ab,ac,ad;q)_n}q^{n^2/2-n/2}(1-abcdq^{2n-1})(-ax)^np_n(y;a,b,c,d|q)\nonumber
.\end{eqnarray}
\end{lz}
\begin{proof} It suffices to set
 $r=s=4, t=1$, and make the specialization of the parameters
\bnm
(a_1,a_2,a_3)&\to& (q, ae^{i\theta},ae^{-i\theta}),\\
(b_1,b_2,b_3,c)&\to& (ab,ac,ad,abcd/q)
\enm
 in (\ref{dlidi}).  As such,  it is clear that the ${}_5\phi_4$ series on the right-hand side
\begin{align*}
&{}_{5}\phi _{4}\left
[\begin{matrix}q^{-n},q,ae^{i\theta},ae^{-i\theta},abcdq^{n-1}\\ ab,ac,ad,q\end{matrix}; q, q\right]\\
&={}_{4}\phi _{3}\left
[\begin{matrix}q^{-n},ae^{i\theta},ae^{-i\theta},abcdq^{n-1}\\ ab,ac,ad\end{matrix}; q, q\right]=\frac{p_n(y;a,b,c,d|q)}{\poq{ab,ac,ad}{n}}a^n,
\end{align*}
while, in view of Eq.(7.5.2) of \cite{10}, $p_n(\cos(\theta);a,b,c,d|q)$ is nothing but the well-known Askey-Wilson polynomial. Thus we have the claimed.
\end{proof}
In the meantime, we may recover Liu's generating function  of the Askey-Wilson polynomials from Theorem \ref{dlmain} or Corollary \ref{tl213} directly.
\begin{lz}[{\rm \cite[cf.][Proposition 1.8]{liunew}}]Let $y=\cos(\theta)$. Then, for $|ax|<1$, it holds
\begin{align}
\sum_{n=0}^{\infty}\frac{(abcd/q,1/x;q)_n }{(q,ab,ac,ad,abcdx;q)_n}&\frac{1-abcdq^{2n-1}}{1-abcdq^{-1}}(ax)^np_n(y;a,b,c,d|q)\nonumber\\
&=\frac{\poq{abcd,ae^{i\theta},ae^{-i\theta},abx,acx,adx}{\infty}}
{\poq{ab,ac,ad,abcdx,axe^{i\theta},axe^{-i\theta}}{\infty}}.\label{wilson-1}\end{align}
\end{lz}
\begin{proof} It suffices to set
 $r=s=3$  in \eqref{eq334} and then make the specialization of the parameters
\bnm
(A_1,A_2,A_3)&\to& (abcd,ae^{i\theta},ae^{-i\theta}),\\
a&\to&abcd/q,\\
(B_1,B_2,B_3)&\to& (ab,ac,ad).
\enm
  As such,  it is clear that the ${}_5\phi_4$ series on the left-hand side
\begin{align*}
&{}_{5}\phi _{4}\left
[\begin{matrix}q^{-n},abcdq^{n-1},abcd,ae^{i\theta},ae^{-i\theta}\\ abcd,ab,ac,ad\end{matrix}; q, q\right]\\
&={}_{4}\phi _{3}\left
[\begin{matrix}q^{-n},abcdq^{n-1},ae^{i\theta},ae^{-i\theta}\\ ab,ac,ad\end{matrix}; q, q\right]=\frac{p_n(y;a,b,c,d|q)}{\poq{ab,ac,ad}{n}}a^n.
\end{align*} Thus \eqref{wilson-1} follows from \eqref{eq334}, as claimed.
\end{proof}
\subsection{Generalized Rogers-Fine identity} As mentioned earlier, Corollary \ref{dlmain-add} generalizes the  famous Rogers-Fine identity which is useful in proving of theta function identities. The interested reader  may consult \cite{fine,roger} for more details.

\begin{lz}[Rogers-Fine identity: \mbox{\cite[p.15, Eq.(14.1)]{fine} or \cite{roger}}]\label{oneoneonne} For $|x|<1$, it holds
\begin{eqnarray}
(1-x)\sum_{n=0}^{\infty}\frac{\poq{c/x}{n}}{\poq{aq}{n}}x^n
=\sum_{n=0}^{\infty}\frac{\poq{c/a,c/x}{n}}{(aq,xq;q)_{n}}(1-cq^{2n})(ax)^nq^{n^2}.\label{rogerfine}
\end{eqnarray}
\end{lz}

\begin{proof} Actually, it is an immediate consequence of (\ref{dlidi}) under the specialization of  the parameters
  $r=s=2,t=1$ and $(a_1,b_1)\to(q,aq)$. In such case, it is clear that
\[{}_{3}\phi _{2}\left
[\begin{matrix}q^{-n},&q,&cq^{n}\\  &aq,&q\end{matrix}; q, q\right]={}_{2}\phi _{1}\left
[\begin{matrix}q^{-n},&cq^{n}\\ &aq\end{matrix}; q, q\right],\]
while we can evaluate the ${}_{2}\phi _{1}$ series on the right-hand side by the $q$-Chu-Vandermonde formula \cite[(II.6)]{10}. This finishes the proof of (\ref{rogerfine}).
 \end{proof}

Apart from this result,  we readily  extend by using  Corollary \ref{dlmain-add} the above Rogers-Fine identity  to the following transformation.

\begin{lz}[Generalized Rogers-Fine identity]\label{tl2} For $|x|<1$, it holds
\begin{eqnarray}
(1-x)\sum_{n=0}^{\infty}\frac{\poq{d,c/x}{n}}{\poq{cd/a,aq}{n}}x^n=\sum_{n=0}^{\infty}\frac{\poq{c/a,aq/d,c/x}{n}}
{\poq{aq,cd/a,xq}{n}}q^{n^2/2-n/2}(1-cq^{2n})(-dx)^n. \label{grfm}
\end{eqnarray}
\end{lz}
\begin{proof} To establish (\ref{grfm}),  we first set
 $r=s=3,t=1$, and make the specialization of the parameters
\bnm
(a_1,a_2)&\to& (q, d),\\
(b_1,b_2)&\to& (aq, cd/a)
\enm
 in (\ref{dlidi}).  Accordingly, it is easy to see that
\[{}_{4}\phi _{3}\left
[\begin{matrix}q^{-n},&q,&d,&cq^{n}\\ &aq,&cd/a,&q\end{matrix}; q, q\right]={}_{3}\phi _{2}\left
[\begin{matrix}q^{-n},&d,&cq^{n}\\ &aq,&cd/a\end{matrix}; q, q\right],\]
while the last ${}_{3}\phi _{2}$ series on the right-hand side  can now be evaluated in closed form by the $q$-Pfaff-Saalsch\"{u}tz ${}_3\phi_2$ summation formula \eqref{saal}.
This gives the complete proof of (\ref{grfm}).
\end{proof}

As is to be expected, upon letting $d$ tend to zero on both sides of (\ref{grfm}), we thereby recover  the Roger-Fine identity (\ref{rogerfine}). As a matter of fact,  by means of Formula \ref{tl2}, we may further deduce a more useful transformation.
\begin{lz}
 Replacing $d$  with $adq/c$ in (\ref{grfm}), we have
 \begin{align}
(1-x)\sum_{n=0}^{\infty}\frac{\poq{adq/c,c/x}{n}}{\poq{dq,aq}{n}}x^n=
\sum_{n=0}^{\infty}\frac{\poq{c/a,c/d,c/x}{n}}
{\poq{aq,dq,xq}{n}}
\left(-\frac{axd}{c}\right)^nq^{n^2/2+n/2}(1-cq^{2n}).\label{grfm-1}
\end{align}
\end{lz}
As might be observed,  the expression on the right-hand side of (\ref{grfm-1}) is symmetric in  $a,d,x$. This property will be used to the proof of Ramanujan's reciprocal theorem below. It should be pointed out that \eqref{grfm-1} is the special case of Watson's transformation from the VWP ${}_8\phi_7$ to  ${}_4\phi_3$ series \cite[(III.18)]{10}, because  (\ref{grfm-1}) is equivalent to
 \begin{align*}
\frac{1-x}{1-c}\,{}_{3}\phi _{2}\left[\begin{matrix}q,&adq/c,&c/x\\ &aq,&dq\end{matrix}
; q, x\right]=\lim_{n\to \infty}\,_{8}W_{7}(c;c/a,c/d,c/x, q,q^{-n};q,dxyq^{1+n}/c).
\end{align*}

\subsection{Andrews' four-parametric  reciprocity theorem and Ramanujan's  ${}_1\psi_1$ summation formula}
As is to be expected,  Formula \ref{tl2} is very helpful in showing  Andrews' four-parametric  reciprocity theorem. As is displayed in \cite[Theorem 1.1]{berndt}, Andrews' four-parametric  reciprocity theorem is a generalization of
 Ramanujan's  reciprocity theorem. The reader may consult  \cite{xrma2} for more details.

\begin{lz}[Andrews' reciprocity theorem: {\rm \cite[Theorem 6]{andrews2}}]\label{1.2} Let $a,d,x,aq/c,$ $dq/c,xq/c$ be any complex numbers other than of the  form $q^{-n}$ for positive integers $n$, $\max\{|x|,|xq/c|\}<1$. Then there holds
 \begin{eqnarray}
\sum_{n=0}^{\infty}\frac{\poq{adq/c,c/x}{n}}{\poq{d,a}{n+1}}x^n&-&\frac{q}{c }\sum_{n=0}^{\infty}\frac{\poq{adq/c,q/x}{n}}{\poq{dq/c,aq/c}{n+1}}\bigg(\frac{xq}{c}\bigg)^n\nonumber\\
&=&\frac{\poq{q,c,q/c,adq/c,dxq/c,axq/c}{\infty}}
{\poq{d,a,x,dq/c,aq/c,xq/c}{\infty}}.\label{kang1}
\end{eqnarray}
 \end{lz}
\begin{proof} At first, by Bailey's VWP $\,_6\psi_6$
summation formula for bilateral series \cite[(II.33)]{10},  we easily obtain
\begin{align*}
\,_6\psi_6\left[%
\begin{array}{ccccccc}
 q\sqrt{c},&-q\sqrt{c}, & c/a, &c/d, & c/x,&y\\
 \sqrt{c},&-\sqrt{c}, & aq, & dq, & xq,& cq/y\\
\end{array};q,\frac{adxq}{cy}\right]\\
=\frac{\poq{q,cq,q/c,adq/c,dxq/c,axq/c,dq/y,aq/y,xq/y}{\infty}}
{\poq{dq,aq,xq,cq/y,dq/c,aq/c,xq/c,q/y, adxq/(cy)}{\infty}}.
 \end{align*}
Define now, for brevity,
 \begin{eqnarray}
H(a,d,x;c):=
\sum_{n=0}^{\infty}\frac{\poq{c/a,c/d,c/x}{n}}
{\poq{aq,dq,xq}{n}}
\left(\frac{axdq}{c}\right)^n\tau(n)(1-cq^{2n}).
\end{eqnarray}
 Recasting  Bailey's VWP ${}_6\psi_6$ summation formula above in terms of $H(a,d,x;c)$, we readily  get
  \begin{eqnarray}
&&\frac{1}{1-c}H(a,d,x;c)-\frac{\kappa qc^2}{c-1}H\bigg(\frac{aq}{c},\frac{dq}{c},\frac{xq}{c};\frac{q^{2}}{c}\bigg)\nonumber\\&=&\lim_{y\to\infty}\,_6\psi_6\left[%
\begin{array}{ccccccc}
 q\sqrt{c},&-q\sqrt{c}, & c/a, &c/d, & c/x,&y\\
 \sqrt{c},&-\sqrt{c}, & aq, & dq, & xq,& cq/y\\
\end{array};q,\frac{adxq}{cy}\right]\nonumber\\
&=&\frac{\poq{q,cq,q/c,adq/c,dxq/c,axq/c}{\infty}}
{\poq{dq,aq,xq,dq/c,aq/c,xq/c}{\infty}},\label{reciprocal}
 \end{eqnarray}
where $\kappa=(a-1) (d-1) (x-1)/((c-a q) (c-d
   q) (c-q x))$. Referring to (\ref{grfm-1}), it is immediate that \begin{eqnarray*}
 H(a,d,x;c)=(1-x)\sum_{n=0}^{\infty}\frac{\poq{adq/c,c/x}{n}}{\poq{dq,aq}{n}}\,x^n
 \end{eqnarray*}
 and correspondingly
 \begin{eqnarray*}
 H\bigg(\frac{aq}{c},\frac{dq}{c},\frac{xq}{c};\frac{q^{2}}{c}\bigg)=(1-xq/c)\sum_{n=0}^{\infty}\frac{\poq{adq/c,q/x}{n}}{\poq{dq^{2}/c,aq^{2}/c}{n}}(xq/c)^n.
 \end{eqnarray*}
 Substitute these two expressions into the left-hand side of (\ref{reciprocal}).
 After some simplification, we obtain the desired identity.
 \end{proof}

It is noteworthy that the limitation as  $d$ to $0$ in  (\ref{kang1}) of Formula \ref{1.2} leads  to  Ramanujan's well known ${}_1\psi_1$ summation formula  \cite[(II.29)]{10}.
\begin{lz}[Ramanujan's ${}_1\psi_1$ summation formula] For $|c|<|x|<1$, it holds
 \begin{eqnarray}
  \sum_{n=-\infty}^{\infty}\frac{\poq{c/x}{n}}{\poq{aq}{n}}\,x^n=\frac{(q, c,q/c,axq/c;q)_{\infty}}
{(x,aq,aq/c,xq/c;q)_{\infty}}.\label{15}
     \end{eqnarray}
 \end{lz}   \begin{proof}
 To establish (\ref{15}), we first take the limit  of both sides of (\ref{kang1}) as $d$ tends to zero, arriving at
 \begin{align}
\sum_{n=0}^{\infty}\frac{\poq{c/x}{n}}{\poq{a}{n+1}}\,x^n-\frac{q}{ c}\sum_{n=0}^{\infty}\frac{\poq{q/x}{n}}{\poq{aq/c}{n+1}}\,\bigg(\frac{xq}{c}\bigg)^n=\frac{\poq{q,c,q/c,axq/c}{\infty}}
{\poq{a,x,aq/c,xq/c}{\infty}}.\label{triple}
\end{align}
Observe that  the identity (\ref{rogerfine}) in Formula \ref{oneoneonne} is  symmetric in two variables $x$ and $a$,  resulting in
\[\sum_{n=0}^{\infty}\frac{\poq{c/x}{n}}{\poq{a}{n+1}}\,x^n=\sum_{n=0}^{\infty}\frac{\poq{c/a}{n}}{\poq{x}{n+1}}\,a^n.\]
Next, by specializing
  $(a,x,c)\to (aq/c,xq/c,q^{2}/c)$, we have
\[\sum_{n=0}^{\infty}\frac{\poq{q/x}{n}}{\poq{aq/c}{n+1}}\bigg(\frac{xq}{c}\bigg)^n
=\sum_{n=0}^{\infty}\frac{\poq{q/a}{n}}{\poq{xq/c}{n+1}}\bigg(\frac{aq}{c}\bigg)^n.\]
Substitute this into (\ref{triple}) and multiply both sides by $1-a$. We therefore obtain
\begin{eqnarray*}
\sum_{n=0}^{\infty}\frac{\poq{c/x}{n}}{\poq{aq}{n}}x^n+\sum_{n=1}^{\infty}
\frac{\poq{1/a}{n}}{\poq{xq/c}{n}}\bigg(\frac{aq}{c}\bigg)^n
=\frac{\poq{q,c,q/c,axq/c}{\infty}}
{\poq{aq,x,aq/c,xq/c}{\infty}}.
\end{eqnarray*}
Apparently, these two series  on the  left-hand side  can now be unified into a bilateral series
\[\sum_{n=-\infty}^{\infty}\frac{\poq{c/x}{n}}{\poq{aq}{n}}x^n,\]
yielding (\ref{15}). \end{proof}

\begin{remark}
It should be mentioned  that in her derivation of Theorem 4.1 in \cite{kang}, S. Y. Kang pointed out that Ramanujan's ${}_1\psi_1$ summation formula  is equivalent to a three-variable generalization of Ramanujan's reciprocity theorem.
\end{remark}
\subsection{Three transformations for VWP
${}_8\phi_{7}$ series}

Recall that there is a symmetric transformation presented in our latest paper \cite{mawang}, wherein it reveals an equivalency of the Bailey VWP ${}_6\psi_6$ summation formula to the Weierstrass  theta function identity.
As we shall show below, it is a special consequence of  Theorem \ref{dlmain-more}.
\begin{lz}[{\rm\cite[Lemma 5]{mawang}}]
Suppose that the parameters $a,b,c,d,e$ are
such that the denominator factors in the terms of the series are never zero such that $|c|>\max\{|q|,|aq|\}$, $t=a^2/(bde)$. Then it holds
\begin{align}
&\sum_{n=0}^\infty\frac{1-aq^{2n+2}/t}{\poq{aq/c}{n+1}}\frac{\poq{cq/t,bq,dq,eq}{n}}{\poq{aq^2/(tb),aq^2/(td),aq^2/(te)}{n}}
\bigg(\frac{q}{c}\bigg)^n\nonumber\\
=&\sum_{n=0}^\infty
\frac{1-q^{2n+2}/t}{\poq{q/c}{n+1}}\frac{\poq{cq/t,bq/a,dq/a,eq/a}{n}}
{\poq{aq^2/(tb),
aq^2/(td),aq^2/(te)}{n}}\bigg(\frac{aq}{c}\bigg)^n.\label{mainid}
\end{align}
\end{lz}
\begin{proof}   It suffices to specialize Theorem \ref{dlmain-more} to   $$(a,b,c,x)\to (aq,q,q/t,1/c),~~ t=a^2/(bde),$$
and
\[\beta_n=\frac{\poq{bq/a,dq/a,eq/a}{n}}
{\poq{aq^2/(tb),
aq^2/(td),aq^2/(te)}{n}}(1-q^{2n+2}/t)(aq)^n.\]
Then, using Jackson's ${}_8\phi_7$ summation formula \cite[(II.22)]{10} for VWP series, we compute
\begin{align*}
 \alpha_n&=\frac{q^{n}\poq{a}{n}(1-aq^{2n+2}/t)}{\poq{q^3/t}{n}}{}_8W_7(q^2/t;q^{-n},aq^{n+2}/t,bq/a,dq/a,eq/a;q,q)\\
 &=\frac{q^{n}\poq{a}{n}(1-aq^{2n+2}/t)}{\poq{q^3/t}{n}}\frac{\poq{q^3/t,bq,dq,eq}{n}}{\poq{a,aq^2/(tb),aq^2/(td),aq^2/(te)}{n}}\\
 &=\frac{\poq{bq,dq,eq}{n}}{\poq{aq^2/(tb),aq^2/(td),aq^2/(te)}{n}}(1-aq^{2n+2}/t)q^{n}.
\end{align*}
A substitution of  $(\alpha_n, \beta_n)$ into \eqref{dl4115} of Theorem \ref{dlmain-more} yields  \eqref{mainid}.\end{proof}
A combination of such pair $(\alpha_n,\beta_n)$ with the WP Bailey lemma given by Corollary \ref{wpbaileylemma} yields directly
\begin{lz} Let $t=a^2/(bde)$. Then
\begin{align}
&{}_{8}W _{7}\left(aq^2/t;x,y,bq,dq,eq; q, q^3/(txy)\right)\nonumber\\
&=\frac{\poq{aq^3/t,q^3/tx,q^3/ty,aq^3/(txy)}{\infty}}{\poq{q^3/t,q^3/(txy),aq^3/tx,aq^3/ty}{\infty}}\,{}_{8}W _{7}\left(\frac{q^2}{t};x,y,\frac{bq}{a},\frac{dq}{a},\frac{eq}{a};q,\frac{aq^3}{txy}\right).\nonumber
\end{align}
\end{lz}

In much the same way as above,  it is not hard to  combine  Gasper and Rahman's  ${}_8\phi_7$ summation formula (II.16) of \cite{10}, we also obtain a transformation from a VWP ${}_8\phi_7$ series with base $q^2$ to one with base $q$. To the best of our knowledge, it does not fall under the purview of Bailey's well-known ${}_{10}\phi_{9}$ transformation \cite[Eq.(2.9.1)]{10} and its limit case ${}_8\phi_7$ transformation \cite[Eq.(2.10.1)]{10}. Also it seems to have been unknown so far in the literature.
\begin{lz} For $|bx|<1,$ it holds
\begin{align}
{}_{8}W_{7}(c;q,1/b^2,c/x,cq/x,q^2;q^2,b^2x^2)\label{ccc}\\
=\frac{(1-x)(1-bc)}{(1-c)(1-bx)}
\,{}_{8}W_{7}(bc;q,b\sqrt{c},-b\sqrt{c},1/b,c/x;q,-bx)\nonumber.
\end{align}
\end{lz}
\begin{proof}
In (\ref{utrans-I}) of Corollary \ref{tl4}, set  $r=5$  and
make
 the specialization of the parameters
 \bnm
(a_1,a_2,a_3)&\to &(w,-w,bcq/w^2),\\
(a_4,a_5)&\to&(q\sqrt{bc},-q\sqrt{bc}),
\enm
and
\[a_1b_1=a_2b_2=a_3b_3=a_4b_4=a_5b_5=bcq.\]
And further set
$$(t,w)\to (-b,-b\sqrt{cq}).$$
With these choices, we are able to apply Gasper and Rahman's  ${}_8\phi_7$ summation formula (II.16) of \cite{10} to evaluate the series  on the right-hand side of (\ref{utrans-I}) in closed form. In the sequel, we obtain that
 \begin{align*}
 &{}_{8}W_{7}(bc;q^{-n},cq^{n},w,-w, bcq/w^{2};q,-bq)
\\ &=\frac{\poq{bcq,w^{2}/bc}{\infty}\poqq{q^{1-n},cq^{n+1},w^{2}q^{n+1},w^{2}q^{-n+1}/c}{\infty}}{\poq{bcq^{n+1},
bq^{1-n}}{\infty}\poqq{q,cq,w^{2}q,w^{2}q/c}{\infty}},
\nonumber
\end{align*}
which further turns out that for  $n=2m$,
 \begin{align*}
 &{}_{8}W_{7}(bc;q^{-n},cq^{n},w,-w, bcq/w^{2};q,-bq)\\
 &=\frac{\poq{bcq}{2m}}{\poq{bq^{1-2m}}{2m}}
\frac{\poqq{q^{1-2m},b^2q^{2-2m}}{m}}{\poqq{b^2cq^2,cq}{m}}=\frac{\poq{bcq}{2m}}{\poq{1/b}{2m}}
\frac{\poqq{q,1/b^2}{m}}{\poqq{b^2cq^2,cq}{m}}\nonumber
\end{align*}
and $0$ otherwise, because  that $\poqq{q^{1-n}}{\infty}=0$ when $n$ is odd. Substitute this result into (\ref{utrans-I}) and reformulate the resulting identity in terms of standard notation of $q$-series. It gives the complete proof of (\ref{ccc}).
\end{proof}

\subsection{Identities arising from some well known formulas}
As illustrations of  Corollaries \ref{dlmain-old} and \ref{dlmain-add},  we can deduce two new series expansions from Gauss' ${}_2\phi_1$ summation formula and one from Rogers' VWP ${}_6\phi_5$ summation formula.
\begin{lz} For  $x,y\neq 1$, we have
\begin{align}
\frac{\poq{ax,ay}{\infty}}{\poq{a,axy}{\infty}}=\sum_{n_1,n_2\geq 0}\tau(n_1+n_2)q^{-n_1n_2}(1-q^{2n_1})(1-q^{2n_2})x^{n_1}y^{n_2}\nonumber\\
\qquad~\quad\times\,\,\frac{\poq{1/x}{n_1}\poq{1/y}{n_2}}{\poq{x}{n_1+1}\poq{y}{n_2+1}}
{}_4\phi_3\left
[\begin{matrix}q^{-n_1},& q^{n_1},&q^{-n_2},& q^{n_2},\\ &q,&q,&a\end{matrix}; q, aq^2\right].\label{id323}
\end{align}

\end{lz}
\begin{proof} Actually, the identity (\ref{id323}) follows immediately  by applying  Corollary  \ref{dlmain-add} twice to  Gauss' ${}_2\phi_1$ summation formula, which saying
 \begin{align}
 \frac{\poq{ax,ay}{\infty}}{\poq{a,axy}{\infty}}={}_2\phi_1\left
[\begin{matrix}1/x,&1/y\\ &a\end{matrix}; q, axy\right].\label{idgauss}
 \end{align}
 The routine calculation is left to the interested reader.
 \end{proof}

Instead, on applying Theorem \ref{dlmain-more} directly to the right-hand function of \eqref{idgauss}, we obtain
\begin{lz}For any complex number $|b|<\infty$, we have
\begin{align}
\frac{\poq{ax,ay}{\infty}}{\poq{a,axy}{\infty}}
=\sum_{n=0}^\infty{}_3\phi_2\left
[\begin{matrix}q^{-n},& bq^{n},&1/y\\ &q,&a\end{matrix}; q, \frac{aqy}{b}\right]\frac{\poq{1/x}{n}}{\poq{bx}{n+1}}
(1-bq^{2n})\tau(n)(bx)^{n}.\label{id323-323}
\end{align}
\end{lz}
\begin{proof} It suffices to reformulate \eqref{idgauss} in the form
 \begin{align*}
 \frac{\poq{ax,ay}{\infty}}{\poq{a,axy}{\infty}}=\sum_{n=0}^\infty\beta_n\,\poq{1/x}{n}x^n
 \end{align*}
 with
 \[\beta_n=\frac{\poq{1/y}{n}}{\poq{q,a}{n}}(ay)^n.\]
 Next, according to \eqref{useful-1}, we easily find that
  \begin{align*}
 \frac{\poq{ax,ay}{\infty}}{\poq{a,axy}{\infty}}=\sum_{n=0}^\infty\alpha_n \frac{\poq{1/x}{n}x^n}{\poq{bx}{n+1}},
 \end{align*}
 where
 \begin{align*}
 \alpha_n&=\tau(n)(1-bq^{2n})b^n\sum_{i=0}^{n}\frac{\poq{q^{-n},bq^{n}}{i}}{\poq{q}{i}}(q/b)^{i}\frac{\poq{1/y}{i}}{\poq{q,a}{i}}(ay)^i\\
 &=\tau(n)(1-bq^{2n})b^n{}_3\phi_2\left
[\begin{matrix}q^{-n},& bq^{n},&1/y\\ &q,&a\end{matrix}; q, \frac{aqy}{b}\right].
 \end{align*}
 It gives the identity \eqref{id323-323}.
 \end{proof}

 We remark that \eqref{id323-323} may be recognized as generalizations of Gauss' ${}_2\phi_1$ summation formula, since the case   $b\to 0$  reduces to \eqref{idgauss}. In a very similar argument,  we can readily deduce a series expansion from Rogers' VWP ${}_6\phi_5$ summation formula as an illustration of  Corollary \ref{dlmain-old}.
\begin{lz} For $|aqxyz|<1$, we have
\begin{align}
\frac{\poq{aq,aqxy,aqxz,aqyz}{\infty}}{\poq{ax,ay,az,aqxyz}{\infty}}=\sum_{n_1,n_2,n_3\geq 0}(ax)^{n_1}(ay)^{n_2}(az)^{n_3}\prod_{i=1}^3(1-q^{2n_i+1})\nonumber\\
\qquad~\quad\times\frac{\poq{q/a,1/x}{n_1}}{\poq{a,xq}{n_1+1}}
\frac{\poq{q/a,1/y}{n_2}}{\poq{a,yq}{n_2+1}}
\frac{\poq{q/a,1/z}{n_3}}{\poq{a,zq}{n_3+1}}\label{id324}\\
\times{}~_{12}W_{11}(a;a,a,a,q^{-n_1},q^{-n_2},q^{-n_3},q^{n_1+1},q^{n_2+1},q^{n_3+1}; q, aq).\nonumber
\end{align}
\end{lz}
\begin{proof} Actually,  by applying repeatedly \eqref{utrans} of Corollary  \ref{dlmain-old}  to Rogers' VWP ${}_6\phi_5$ summation formula, namely,
\begin{align}
 \frac{\poq{aq,aqxy,aqxz,aqyz}{\infty}}{\poq{aqx,aqy,aqz,aqxyz}{\infty}}={}_6W_5(a;1/x,1/y,1/z;q,aqxyz),\label{rogerid}
 \end{align}
the identity (\ref{id324}) follows at once.
 The routine calculation is left to the interested reader.
 \end{proof}

An interesting case (i.e., $r=s=2$) covered by  Theorem  \ref{dlmain} is   Jackson's VWP ${}_8\phi_7$ summation formula.  It  utilizes the VWP ${}_8\phi_7$ transformation formula (III.23)  of \cite[]{10}.
\begin{lz}[Jackson's VWP ${}_8\phi_7$ summation formula: {\rm \cite[cf.][(II.22)]{10}}]\label{jformula} Let $\lambda=a^2q/(bcd)$. For integers $n$, we have
\begin{align}
{}_{8}W _{7}\left(\lambda;\frac{ \lambda\,b}{a},\frac{ \lambda\,c}{a},\frac{ \lambda\,d}{a},aq^n,q^{-n};q,q\right)=\frac{\poq{b,c,d,\lambda q}{n}}{\poq{aq/b,aq/c,aq/d,a/\lambda}{n}}.\nonumber
\end{align}
\end{lz}
\begin{proof} Actually, (III.23)  of \cite{10}  states  that for $|\lambda\,q/(ef)|<1$ and $|aq/(ef)|<1,$
\begin{align}
&{}_{8}W _{7}\left(a;b,c,d,e,f; q,\lambda\,q/(ef)\right)\nonumber\\
&=\frac{\poq{a q, a q/(e f), \lambda\,q/e, \lambda\,q/f}{\infty}}{\poq{aq/e, aq/f, \lambda\,q, \lambda\,q/(ef)}{\infty}}\,{}_{8}W _{7}\left(\lambda;\frac{ \lambda\,b}{a},\frac{ \lambda\,c}{a},\frac{ \lambda\,d}{a},e,f;q,\frac{aq}{ef}\right).\label{2.6}
\end{align}
Once taking $f=1/x$ and defining
\[K(x):=\frac{\poq{ aq,\lambda\,q/e,\lambda\,qx,a qx/e}{\infty}}{\poq{\lambda\,q,aq/e,aqx, \lambda\,qx/e}{\infty}},\]
we can restate \eqref{2.6} in  notation of  Theorem \ref{dlmain} as the following expansion formula:
\begin{align}
F(x)=\sum_{n=0}^\infty \alpha_n \frac{\poq{1/x}{n}}{\poq{aqx}{n}}x^n,\label{4.3}
\end{align}
where
 \begin{align}
F(x)=K(x)\sum_{n=0}^\infty \beta_n \frac{\poq{1/x}{n}}{\poq{\lambda\,qx}{n}}x^n. \label{4.33}
\end{align}
It is very clear from \eqref{2.6} that the coefficients of the expansion \eqref{4.3} and \eqref{4.33} ought to be
\begin{align*}
 \alpha_n&=\frac{1-aq^{2n}}{1-a}\frac{\poq{a,b,c,d,e}{n}}{\poq{q,aq/b,aq/c,aq/d,aq/e}{n}}\bigg(\frac{\lambda\,q}{e}\bigg)^{n},\\
\beta_{n}&=\frac{1-\lambda\,q^{2n}}{1-\lambda}\frac{\poq{\lambda,\lambda\,b/a,\lambda\,c/a,\lambda\,d/a,e}{n}}{\poq{q,aq/b,aq/c,aq/d,\lambda\,q/e}{n}}\bigg(\frac{aq}{e}\bigg)^{n}.
\end{align*}
On the other hand, by applying Lemma \ref{yl3} to \eqref{4.3}, we can evaluate
directly \begin{align*}
 \alpha_n
&=\frac{1-aq^{2n}}{1-a}\frac{(a;q)_{n}}{\poq{q}{n}}
\sum_{k=0}^nq^k
\frac{\poq{q^{-n},aq^n}{k}}{\poq{q,aq}{k}}K(q^k)
\sum_{i=0}^k\beta_i\frac{\poq{q^{-k}}{i}}{\poq{\lambda\,q^{k+1}}{i}}
q^{ki}.
\end{align*}
Note that
\[K(q^k)=\frac{\poq{aq,\lambda\,q/e}{k}}{\poq{\lambda\,q,aq/e}{k}}.\]
We get
\begin{align}
 \alpha_n
&=\frac{1-aq^{2n}}{1-a}\frac{(a;q)_{n}}{\poq{q}{n}}\sum_{i=0}^n\beta_iq^i\tau(i)\frac{\poq{q^{-n},aq^n,\lambda\,q/e}{i}}{\poq{aq/e}{i}\poq{\lambda\,q}{2i}}W_{n,i},\label{finalresult}
\end{align}
where  the inner sum, which can be evaluated by the $q$-Pfaff-Saalsch\"{u}tz ${}_3\phi_2$ summation formula \eqref{saal}, turns out to be
$$
W_{n,i}:={}_{3}\phi _{2}\left[\begin{matrix}q^{-n+i},aq^{n+i},\lambda\,q^{i+1}/e\\  \lambda\,q^{2i+1},aq^{i+1}/e\end{matrix}
; q, q\right]=\frac{\poq{eq^i,q^{i-n+1}\lambda/a}{n-i}}{\poq{\lambda\,q^{2i+1},q^{-n}e/a}{n-i}}.
$$
A direct substitution of this expression as well as those for $\alpha_n$ and $\beta_n$ into \eqref{finalresult} gives the desired identity.
 \end{proof}
\subsection{Further curious $q$-series identities}
As is expected to be, it looks quite promising  to seek for new  $q$-series identities  by applying our expansion formulas to certain existing
 summation formula. Due to space limitations, we shall only consider some typical cases.
 \begin{lz}\label{tl00} Assume that all sums converge. Then we have
\begin{align}
\frac{\poqq{bx,xq/b}{\infty}}{\poq{x/q}{\infty}}=\sum_{n=0}^{\infty}{}_{3}\phi _{3}\left
[\begin{matrix}q^{-n},&b,&q/b\\ q^{-n}/a,&q,&-q\end{matrix}; q, -\frac{q}{a}\right]\frac{(aq;q)_{n} }{(ax;q)_{n+1}}(x/q)^{n}\label{added1},\\
\frac{\poqq{bx,xq/b}{\infty}}{\poq{x/q}{\infty}}=\sum_{n=0}^{\infty}{}_{3}\phi _{2}\left
[\begin{matrix}q^{-n},&b,&q/b\\ &q,&-q\end{matrix}; q, -q^{n+1}\right](x/q)^{n}.\label{added0}
\end{align}
Furthermore, for $b\neq q$,
\begin{align}
\frac{\poqq{q/b}{n}}{\poqq{q^2}{n}}\,{}_{2}\phi _{1}\left
[\begin{matrix}q^{-2n},&bq\\ &bq^{1-2n}\end{matrix}; q^2, b\right]=q^{-n}\,{}_{3}\phi _{2}\left
[\begin{matrix}q^{-n},&b,&q/b\\ &q,&-q\end{matrix}; q, -q^{n+1}\right].\label{added00}
\end{align}
\end{lz}
\begin{proof} To achieve (\ref{added1}), it suffices to set
 $r=2,s=3$ and $t=-1$ in (\ref{utrans}) and then make  the following specialization of the parameters
\bnm
(a_1,a_2)&\to &(b,q/b),\\
(b_1,b_2,b_3)&\to&(q,-q,0),\\
(b,c)&\to&(1/q,0).
\enm
Upon making use of the $q$-analogue of Bailey's ${}_2F_1(-1)$ summation formula (II.10) of \cite{10}, we have
\begin{eqnarray}
{}_{3}\widetilde{\phi}_{4}\left
[\begin{matrix}b,&q/b,&0\\ 0,\,\, q,&-q,&x\end{matrix} ; q, -x\right]={}_{2}\phi _{2}\left
[\begin{matrix}b,&q/b\\ -q,&x\end{matrix}; q, -x\right]=\frac{\poqq{bx,xq/b}{\infty}}{\poq{x}{\infty}},\nonumber\end{eqnarray}
reducing (\ref{utrans})  to (\ref{added1}). Start with (\ref{added1}) and let $a$ tend to zero on both sides. Then (\ref{added0}) is obtained. By equating the coefficients of $x^n$ on both sides of (\ref{added0}), we further obtain (\ref{added00}).
\end{proof}

By combining Corollary \ref{dlmain-add}  with Andrews' terminating summation formula \cite[(II.17)]{10},  we immediately obtain
\begin{lz}\label{tl3}For $|x|<1$, we have
\begin{eqnarray}
&&(1-x)\sum_{n=0}^{\infty}
\frac{\poq{a,-a,c^2q/x}{n}}{\poq{cq,-cq,a^2}{n}}\,x^n\label{tl3id}\\
&=&\sum_{n=0}^{\infty}\frac{(q,c^2q^2/a^2,c^2q/x,c^2q^2/x;q^2)_{n}}
{(c^2q^2,a^2q,xq,xq^2;q^2)_{n}}(1-c^2q^{4n+1})q^{2n^2-n}(ax)^{2n}.\nonumber
\end{eqnarray}
\end{lz}
\begin{proof} To show (\ref{tl3id}), we only need to set
 $r=s=4$, and $t=1$ in (\ref{dlidi}) and then make the simultaneous  specialization of the parameters
\bnm
(a_1,a_2,a_3)&\to &(q,a,-a),\\
 c&\to & c^2q,\\
(b_1,b_2,b_3)&\to&(cq,-cq,a^2).
\enm
In such case, it readily follows  that
\begin{eqnarray*}
&&{}_{5}\phi _{4}\left
[\begin{matrix}q^{-n},&q,&a,&-a,&c^2q^{n+1}\\ &cq,&-cq,&a^2,&q\end{matrix}; q, q\right]\\
&=&{}_{4}\phi _{3}\left
[\begin{matrix}q^{-n},&a,&-a,&c^2q^{n+1}\\ &cq,&-cq,&a^2\end{matrix}; q, q\right]
=\left\{\displaystyle
     \begin{array}{ll}
       0, \,\,&\mbox{if}\,\, n=2m+1 \\
&\\
    \displaystyle   \frac{a^{2m}(q,c^2q^2/a^2;q^2)_{m}}{(c^2q^2, a^2q;q^2)_{m}},  &\mbox{if}\,\,n=2m.
     \end{array}
   \right.\nonumber
\end{eqnarray*}
Note that the last equality is given by Andrews' terminating summation formula (II.17) of \cite{10}.
 All these together gives rise to (\ref{tl3id}).
\end{proof}

Using the $q$-analogue of Whipple's ${}_3F_2$ summation formula \cite[(II.19)]{10}, we may easily find that
\begin{lz}For $|x|<1$, we have
\begin{eqnarray}
&&\sum_{n=0}^{\infty}\frac{\poq{a,-a,q/x}{n}}
{(-q,e,a^2q/e;q)_{n}}
\,x^n\label{tl345}\\
&=&\sum_{n=0}^{\infty}\frac{(eq^{-n},eq^{n+1},a^2q^{1-n}/e,a^2q^{n+2}/e;q^2)_{\infty}}
{(e,a^2q/e;q)_{\infty}}
\frac{(q/x;q)_{n} }{(x;q)_{n+1}}~(q^{n^2}-q^{(n+1)^2})(-x)^{n}.\nonumber
\end{eqnarray}
\end{lz}
\begin{proof} It suffices to set
 $r=s=4$ and $t=1, c=q$ in (\ref{dlidi}) and specialize  the parameters simultaneously
 \bnm
(a_1,a_2,a_3)&\to &(q,a,-a),\\
(b_1,b_2,b_3)&\to&(e,a^2q/e,-q).
\enm
In this case, the ${}_{5}\phi _{4}$ series in (\ref{dlidi})
\[{}_{5}\phi _{4}\left
[\begin{matrix}q^{-n},&q,&a,&-a,&q^{n+1}\\ &e,&a^2q/e,&-q,&q\end{matrix}; q, q\right]={}_{4}\phi _{3}\left
[\begin{matrix}q^{-n},&a,&-a,&q^{n+1}\\ &e,&a^2q/e,&-q\end{matrix}; q, q\right]\]
while the ${}_{4}\phi _{3}$ series can be evaluated by (II.19) of \cite{10} to be
\[{}_{4}\phi _{3}\left
[\begin{matrix}q^{-n},&a,&-a,&q^{n+1}\\ &e,&a^2q/e,&-q\end{matrix}; q, q\right]=\frac{(eq^{-n},eq^{n+1},a^2q^{1-n}/e,a^2q^{n+2}/e;q^2)_{\infty}q^{\binom{n+1}{2}}}
{(e,a^2q/e;q)_{\infty}}.\]
A direct substitution of this into (\ref{dlidi}) yields (\ref{tl345}).
\end{proof}

Continuing along this line, we  readily deduce two concrete identities from Corollary \ref{tl4} via the use of  other known summation  formulas.
\begin{lz} For $\max\{|x|,|cx/q|\}<1$, it holds
\begin{eqnarray}
\sum_{n=0}^\infty \frac{\poq{c/x}{n}}{\poq{cx/q}{n+1}}x^n\label{utrans-useful}=
\sum_{n=0}^{\infty}
\frac{\poq{q^{n+1},c,c/x}{n}}
{\poq{c^2/q}{2n+1}(x;q)_{n+1}}(1-cq^{2n})(cx/q)^{n}.
\end{eqnarray}
\end{lz}
\begin{proof}
It suffices to set   $r=0, t=1,$ and
$b=c/q$ in (\ref{utrans-I}) of Corollary \ref{tl4}. Then the $q$-Pfaff-Saalsch\"{u}tz ${}_3\phi_2$ summation formula \eqref{saal} is applicable to the ${}_{3}\phi _{2}$ series on the right-hand side of \eqref{utrans-I}, yielding
\begin{eqnarray}
{}_{3}\phi _{2}\left
[\begin{matrix}q^{-n},&c^2/q,&cq^{n}\\ &c^2q^{n},&cq^{-n}\end{matrix}; q, q\right]=\frac{\poq{c,q^{n+1}}{n}}
{\poq{c^2q^{n},q/c}{n}}.\nonumber\end{eqnarray}
Thus we get the desired transformation.
\end{proof}

It is of interest to note that the limit as $x$ tends to zero offers a new identity for the partial theta function. Regarding the latter, the reader may consult \cite{andrews2,warnaar} for further information.

\begin{lz} We have
\begin{eqnarray}
\sum_{n=0}^\infty (-1)^nq^{n^2/2-n/2}y^n
=\sum_{n=0}^{\infty}
\frac{\poq{q^{n+1},y}{n}}
{\poq{y^2/q}{2n+1}}(-1)^nq^{n^2/2-3n/2}(1-yq^{2n})y^{2n}.\label{12345}
\end{eqnarray}
\end{lz}
\begin{proof} At first, taking the limit $x\to 0$ of (\ref{utrans-useful}) yields
\begin{eqnarray*}
\sum_{n=0}^\infty \tau(n)c^n=\sum_{n=0}^{\infty}\tau(n)(1-cq^{2n})(c^2/q)^{n}
\frac{\poq{q^{n+1},c}{n}}
{\poq{c^2/q}{2n+1}}.
\end{eqnarray*}
And then replacing $c$ with $y$, we get
(\ref{12345}).
\end{proof}

\section{Concluding remarks}
\setcounter{equation}{0}
Thus far we have only exploited applications of the $(1-xy,y-x)$-expansion formula with  appropriate
parametric specializations to  summations and transformations of basic hypergeometric series, wherein (as one of main reasons for us to do so) the absolute and uniform convergence of the  infinite series involved is  easy to justify, thereby  the limit function  is analytic. Nevertheless, as indicated by Remark \ref{rmk},
 this idea is applicable to arbitrary analytic functions with proper  specialization of parameters. For instance, consider
 \begin{equation}\label{00}
    F(x)=\frac{\poq{x}{\infty}}{\poq{ax}{\infty}}{}_3\phi_2\left
[\begin{matrix}a,&b,&c/x\\ &d,&abc/d\end{matrix}; q, x\right]
 \end{equation}
 with the choices that  $x_n=aq^{n-1}$ and $b_n=cq^n$,
 then the $(1-xy,y-x)$-expansion formula will lead us to the limitation of Watson's terminating ${}_{8}\phi_{7}$ to ${}_{4}\phi_{3}$ series \cite[cf.][(III.18)]{10}:
   \begin{eqnarray}
&&\frac{\poq{x,ac}{\infty}}{\poq{ax,c}{\infty}}\,{}_3\phi_2\left
[\begin{matrix}a,&b,&c/x\\ &d,&abc/d\end{matrix}; q, x\right]\\
&=&\lim_{n\to \infty}{}_{8}W_{7}(ac/q;a,d/b,ac/d,c/x,bq^{-n};q,xq^n).\nonumber
\end{eqnarray}
From this  perspective it is reasonable to believe that, under other specialization of parameters just as in (\ref{gessel1})/(\ref{gessel2}), $x_n=ap^n, b_n=cq^n, p\neq q$,  the $(1-xy,y-x)$-expansion formula will result in some transformation formulas for bibasic series \cite[cf.][]{8,rahman}, including  other $(f,g)$-expansion formulas such as Warnaar's matrix inversion \cite[cf.][Lemma 3.2]{1000} and \cite[cf.][]{29} for elliptic hypergeometric series.

Another long-standing problem puzzled the second author may be written like this: does there exist any general (not merely for  special conditions) summation formula for the VWP ${}_{r+1}\phi_r$ series when $r>8$? Our query is based on a phenomenon that in the world of $q$-series, we too often meet and utilize various summation formulas like the ${}_{1}\phi_0$($q$-binomial theorem),
${}_{2}\phi_1$(Chu-Vandermonde, Gauss), ${}_{3}\phi_2$(Pfaff-Saalsch\"{u}tz),  ${}_{4}\phi_3$(Andrews), ${}_{8}\phi_7$(Jackson) and its limit ${}_{6}\phi_5$. In this sense,  Jackson's ${}_{8}\phi_7$ summation formula is  certainly one of the most general (highest \emph{level} so far) formulas in closed form. As a first step toward deep investigation of this problem, all expansions obtained in our foregoing discussion  remind us that if there did exist a closed formula for the ${}_{r+1}\phi_r$ series, then we should have the corresponding one for the ${}_{r}\phi_{r-1}$ series.  All these problems deserve further research.
%\end{document}
\section*{Acknowledgements}
 This  work was supported by the National Natural Science Foundation of China [Grant No.  11471237].

%\section*{References}

\bibliography{wangjin29191222}

\end{document}